\documentclass[11pt,reqno]{amsart}
\usepackage[margin=1in]{geometry}
\usepackage[T1]{fontenc}
\usepackage{lmodern}
\usepackage{microtype}
\usepackage{amsmath,amssymb,amsthm,mathtools}
\usepackage{enumitem}
\usepackage{needspace}
\usepackage{booktabs,array}
\usepackage{tikz}
\usetikzlibrary{arrows.meta,calc,patterns,decorations.pathreplacing}
\usepackage[section]{placeins}
\definecolor{figink}{RGB}{35,42,50}
\definecolor{figmuted}{RGB}{105,110,118}
\definecolor{figfirst}{RGB}{43,88,132}
\definecolor{figsecond}{RGB}{146,95,32}
\definecolor{figcontact}{RGB}{164,49,60}
\definecolor{figactive}{RGB}{27,111,98}
\tikzset{
 figaxis/.style={draw=figmuted,line width=.45pt,-{Stealth[length=1.5mm]}},
 figguide/.style={draw=figmuted!55,densely dashed,line width=.4pt},
 figleader/.style={draw=figmuted,line width=.4pt},
 figlabel/.style={fill=white,inner sep=1.5pt},
 figtitle/.style={font=\small\bfseries,text=figink},
 figtext/.style={font=\footnotesize,text=figink},
 figcurve/.style={draw=figfirst,line width=1.1pt},
 figbound/.style={draw=figink!65,dash pattern=on 3pt off 2pt,line width=.75pt}
}
\usepackage[hidelinks]{hyperref}
\numberwithin{equation}{section}
\newtheorem{maintheorem}{Theorem}

\newtheorem{theorem}{Theorem}[section]
\newtheorem{proposition}[theorem]{Proposition}
\newtheorem{lemma}[theorem]{Lemma}

\theoremstyle{definition}
\newtheorem{definition}[theorem]{Definition}
\newtheorem{remark}[theorem]{Remark}
\newcommand{\C}{\mathbb C}
\newcommand{\R}{\mathbb R}
\newcommand{\PP}{\mathbb P}
\newcommand{\PSH}{\operatorname{PSH}}
\newcommand{\Sch}{\mathcal S}
\newcommand{\diag}{\operatorname{diag}}
\newcommand{\rank}{\operatorname{rank}}

\newcommand{\eps}{\varepsilon}
\newcommand{\omegaFS}{\omega_{\mathrm{FS}}}

\newcommand{\hhat}[1]{\hat{\hat{#1}}}
\title[Spatial degeneracy of weak geodesics]{Spatial degeneracy along weak geodesics in the space of K\"ahler potentials}
\author{Xia Xiao}
\address{Département de mathématiques, Université du Québec à Montréal, Montréal, Canada}
\email{xia.xiao.math@gmail.com}
\author{Qi Yao}
\address{IASM, Zhejiang University, Hangzhou, Zhejiang, China}
\email{qiyao11@zju.edu.cn}
\date{}
\hypersetup{pdfauthor={Xia Xiao and Qi Yao},pdftitle={Spatial degeneracy along weak geodesics in the space of K\"ahler potentials}}
\begin{document}
\raggedbottom
\begin{abstract}
We study spatial Hessian degeneracy along weak K\"ahler geodesics.
We construct examples in three settings:
\begin{enumerate}[label=\arabic*.,leftmargin=*,nosep]
\item On every compact K\"ahler manifold, the spatial form degenerates
on a nonempty open set during an interior time interval, although
both endpoints are smooth and strictly K\"ahler.
\item On ALE K\"ahler manifolds with a holomorphic quotient end and a
decaying potential asymptotic to the Euclidean metric, spatial
degeneracy occurs on open regions arbitrarily far along an unbounded horn.
\item For cone and Poincar\'e endpoint metrics along a smooth divisor,
normal spatial degeneracy occurs on open sets approaching the divisor.
Both cone endpoints retain the prescribed angle.
\end{enumerate}
The first construction answers negatively the spatial-positivity
question recorded in Darvas's survey
\cite[Question~2.3]{Darvas2026}. The ALE and divisor constructions
show that weak geodesics need not preserve the asymptotic behavior
of the endpoint metrics.
\end{abstract}
\maketitle
\markboth{\shorttitle}{\shorttitle}

\section{Introduction and statement of main results}\label{sec:introduction}

Let $(X,\omega)$ be a K\"ahler manifold. Fixing the reference form $\omega$,
we write the space of K\"ahler potentials as follows:
\[
\mathcal H(X,\omega)
=\{u\in C^\infty(X,\R):\omega+i\partial_X\bar\partial_X u>0\}.
\]
In the compact setting, the Riemannian metric on this space was introduced
by Mabuchi \cite[Section 2, equation (2.1.1)]{Mabuchi1987}.
The homogeneous complex Monge--Amp\`ere formulation is due to Semmes
\cite{Semmes1992} and Donaldson
\cite[Section 5, Proposition 3]{Donaldson1999}. More precisely, let
\[
\Sigma=\{\tau\in\C:1<|\tau|<e\},\qquad t=\log|\tau|,\qquad
\Omega=\pi_X^*\omega.
\]
For a path $u_t$, set $U(x,\tau)=u_{\log|\tau|}(x)$. The geodesic equation is the Dirichlet problem of the following HCMA equation:
\begin{equation}\label{eq:hcma}
\begin{cases}
(\Omega+i\partial\bar\partial U)^{n+1}=0&\text{on }X\times\Sigma,\\
\Omega+i\partial\bar\partial U\ge0&\text{on }X\times\Sigma,\\
U=u_j&\text{on }X\times\{|\tau|=e^j\},\quad j=0,1.
\end{cases}
\end{equation}
Here $n=\dim_{\C}X$. The differential operators $\partial$, $\bar \partial$ in \eqref{eq:hcma} act on the product,
whereas $\partial_X$, $\bar\partial_X$ will denote spatial derivatives.

On a compact K\"ahler manifold, Chen \cite[Theorem 3]{Chen2000} proved the existence
of weak geodesics with bounded complex Hessian. Further estimates were
obtained by B{\l}ocki \cite[Theorem 1.4]{Blocki2012}, and full real $C^{1,1}$ regularity
was established by Chu--Tosatti--Weinkove
\cite[Theorem 1.1]{ChuTosattiWeinkove2017}. These estimates do not imply a positive
lower bound for the spatial form $\omega+i\partial_X\bar\partial_X u_t$.
The regularity obstructions are different: Lempert--Vivas
\cite[Theorem 1.2]{LempertVivas2013} exclude $C^3$ solutions of the
constant-rank geodesic boundary problem under their symmetry hypothesis,
and Darvas--Lempert
\cite[Theorem 1.1]{DarvasLempert2012} exclude solutions with continuous
complex Hessian in a related setting. Darvas
\cite[Theorem 1.2]{Darvas2014Morse} gives a more general obstruction to
smooth geodesics. A nonsmooth potential need not have a zero spatial
eigenvalue, and a smooth weak potential may have a degenerate spatial form.

On the initial-value side, Xiao \cite[Theorem B]{Xiao2026}
constructed distinct bounded weak geodesic rays on $\PP^1$
issuing from zero with the same pointwise right initial tangent.

The full space-time determinant in \eqref{eq:hcma} vanishes by the
equation itself. At a smooth point, strict positivity of the spatial
form implies that the full space-time form has rank $n$. Conversely,
for a semipositive form of rank $n$ on the $(n+1)$-dimensional product,
the spatial restriction is positive precisely when its kernel contains
no nonzero spatial tangent vector. This is a linear-algebra condition,
not a consequence of the rank alone.
Our examples have an open region where the spatial rank drops and
the potential is independent of time, so the full rank drops there
as well.

The maximum-rank problem has a long history. Guan--Phong
\cite[Introduction and Theorems 1--2]{GuanPhong2012} proved related
maximum-rank results for real equations and explicitly raised the
complex geodesic question. Ross--Witt Nystr\"om
\cite[Theorem 1.1 and Section 2]{RossWittNystrom2019Rank} constructed
smooth strictly K\"ahler boundary data on $\PP^1\times\partial\Delta$
for which the HCMA solution is completely degenerate on an open subset.
Their parameter domain is a disc. Their boundary family depends on
its angular parameter; the rotation used to relate it to an invariant
problem also introduces a singularity at the centre of the parameter
disc. That construction does not directly give an annular geodesic
segment with two smooth strictly K\"ahler endpoints.

For segments, the spatial-positivity question is recorded explicitly as
a folklore question in \cite[Question 2.3]{Darvas2026}: given smooth
strictly K\"ahler endpoints on a compact manifold, must there be
$\eps>0$ such that
\begin{equation}\label{eq:darvas-question}
\omega+i\partial_X\bar\partial_X u_t\ge\eps\omega,\qquad 0\le t\le1?
\end{equation}
Positive results are known for invariant data on compact smooth toric
varieties \cite[Theorem 3]{Guan1999}, for sufficiently small global
$C^{4,\alpha}$ perturbations \cite[Theorem 1.2]{ChenFeldmanHu2020},
and for sufficiently small global $C^2$ endpoint potentials
\cite[Theorem 1.2]{Hu2024LowerBound}. On compact complex affine manifolds
with a constant-coefficient K\"ahler background, Hu
\cite[Theorem 1.1]{Hu2025Convexity} proves preservation of a common
strict $(S,\omega_0)$-convexity condition, for a fixed constant tensor
$S$. The compact construction below produces a zero spatial eigenvalue
on every compact K\"ahler manifold, with endpoint potentials supported
in any prescribed nonempty open set. Localization of the support does
not impose the global second-order smallness required by the positive
results above.

A related question arises on complete noncompact manifolds. If the
endpoint potentials decay rapidly at infinity, does the spatial form
along the geodesic remain nondegenerate outside a compact set?
Earlier ALE results include Aleyasin
\cite{AleyasinALE} and the second named author
\cite[Theorems A--B]{YaoALEGeodesics};
see also \cite[Theorem 1]{AuvrayPoincare} for Poincar\'e-type geometry.
The weighted estimates for $\varepsilon$-geodesics in
\cite[Theorem B]{YaoALEGeodesics} allow their constants to depend on
$\varepsilon^{-1}$; they do not supply an $\varepsilon$-uniform
second-order asymptotic estimate for the homogeneous limit.
The disc problem has stronger exterior regularity under its stated
hypotheses \cite[Theorem B]{YaoALE}.
We show that, for the annular problem, even Schwartz endpoint decay
does not imply spatial nondegeneracy at infinity.

\subsection{Degeneracy on compact K\"ahler manifolds}

We first construct weak geodesics on compact K\"ahler manifolds for
which the spatial form $\omega+i\partial_X\bar\partial_Xu_t$
degenerates on an open set.

\begin{maintheorem}[Compact K\"ahler manifolds]\label{thm:compact}
Let $(X,\omega)$ be a compact K\"ahler manifold of complex dimension
$n\ge1$. There are smooth real potentials $u_0,u_1$ on $X$, with
$\omega+i\partial_X\bar\partial_Xu_j>0$, having the following
property. For every $0<\delta<1/2$ there is a nonempty open set
$U_\delta\subset X$ such that their weak geodesic is smooth and
independent of time on an open neighborhood of
$U_\delta\times[\delta,1-\delta]$, where
\[
 \rank\bigl(\omega+i\partial_X\bar\partial_Xu_t\bigr)=n-1.
\]
The same endpoint pair works for every $\delta$.

Moreover, for any prescribed nonempty open set $O\subset X$, the
endpoints may be chosen nonnegative with compact support in $O$,
and $U_\delta\Subset O$ for every $\delta$.
\end{maintheorem}

The endpoint and geodesic forms lie in the same prescribed K\"ahler
class $[\omega]$. When $n=1$, the spatial form vanishes on the
indicated region. In every dimension, the theorem answers the
spatial-positivity question~\eqref{eq:darvas-question} negatively on
each compact K\"ahler manifold.
The geodesic has global real $C^{1,1}$ regularity by
\cite[Theorem 1.1]{ChuTosattiWeinkove2017}.
Section~\ref{sec:compact} first constructs model potentials on
$\PP^1$ and proves the corresponding result on K\"ahler products.
It then inserts the same contact configuration into a coordinate
ball to prove Theorem~\ref{thm:compact}.

\subsection{Degeneracy at infinity on ALE spaces}

Let $(X,\omega)$ be a complete K\"ahler manifold of complex dimension
$n\ge1$. Assume that, outside a compact set, $X$ is biholomorphic to
$(\C^n\setminus\overline B_{R_0})/\Gamma$, where
$\Gamma\subset U(n)$ is finite and acts freely on the unit sphere.
Writing $r=|z|$ and $\omega_0=i\partial_X\bar\partial_X|z|^2$, suppose
that on this covering space
\begin{equation}\label{eq:ale-potential}
 \omega=\omega_0+i\partial_X\bar\partial_X\psi,\qquad
 D^k\psi=O(r^{-\sigma-k})\quad(k\ge0),\qquad \sigma>0,
\end{equation}
where $\psi$ is smooth, real and $\Gamma$-invariant. Derivatives are
taken in the holomorphic coordinates of the end.

\begin{maintheorem}[ALE ends]\label{thm:ale}
Under these assumptions, there are smooth bounded real endpoint
potentials $u_0,u_1$ and constants $0<\lambda\le\Lambda$ such that
\[
 \lambda\omega\le\omega+i\partial_X\bar\partial_Xu_j
 \le\Lambda\omega,\qquad
 D^ku_j=O(r^{-\sigma-k})\quad(k\ge0),\qquad j=0,1.
\]
For every compact set $K\subset X$, there are a nonempty open set
$V_K\subset X\setminus K$ and an open interval
$\tfrac12\in I_K\Subset(0,1)$ on whose product their weak geodesic
is smooth, independent of time, and satisfies
\[
 \rank\bigl(\omega+i\partial_X\bar\partial_Xu_t\bigr)=n-1.
\]
The same endpoint pair works for every $K$.
\end{maintheorem}

The construction gives, outside a compact set,
\begin{equation}\label{eq:ale-endpoint-decay}
 u_j=-\psi+s_j,
 \qquad D^ks_j=O(r^{-N})\quad\text{for every }k\ge0\text{ and }N>0.
\end{equation}
Thus the endpoint difference decays rapidly, and both endpoint metrics
have the same leading Euclidean metric as $\omega$.
For $X=\C^n$ and $\omega=\omega_0$, the endpoints are Schwartz functions
and the geodesic is globally $C^{1,1}$. On the open sets in the
construction it is given explicitly by $u_t(z)=-2(\operatorname{Im}z_1)^2$.
Consequently, even Schwartz endpoint decay does not ensure convergence
of the spatial Hessian to the Euclidean one at infinity.

Subsection~\ref{sec:ale} proves the theorem by changing the reference
potential and applying the Euclidean construction on the end.

\subsection{Degeneracy near cone and Poincar\'e divisors}

We next construct rank loss near a smooth divisor. For cone endpoints,
it occurs in the normal direction near a divisor point for a fixed
time interval, while both endpoints retain the prescribed cone angle.
For Poincar\'e endpoints, normal rank loss occurs on open sets
approaching the divisor.

For simplicity, the construction below treats only smooth divisors.
Let $(X,\omega)$ be a compact K\"ahler manifold of complex dimension
$n\ge1$, let $E\subset X$ be a nonempty smooth divisor, and fix
$0<\beta<1$. Following \cite[Section 2, equation (2.1)]{ZhengCone2017}, a K\"ahler cone
metric of angle $2\pi\beta$ along $E$ is a closed positive current with
bounded continuous local potentials, smooth and positive on $X\setminus E$,
and locally comparable in holomorphic coordinates $(z,w)$ with
$E=\{z=0\}$ to
\begin{equation}\label{eq:cone-model}
\omega_{\rm cone}
=\beta^2|z|^{2\beta-2}i\,dz\wedge d\bar z
+i\sum_{\alpha=1}^{n-1}dw_\alpha\wedge d\bar w_\alpha.
\end{equation}
The model has continuous potential $|z|^{2\beta}+|w|^2$, although the coefficient of $dz\wedge d\bar z$ diverges at $E$. We consider the weak geodesic on $X$ relative to the smooth form $\omega$.
The endpoint forms are positive cone metrics; the intermediate forms
below have a vanishing normal coefficient near the chosen divisor point.

\Needspace{6\baselineskip}
\begin{maintheorem}[Cone divisors]\label{thm:conic-boundary}
With $(X,\omega)$, $E$ and $\beta$ as above, fix
$x_*\in E$, an open neighborhood $N$ of $x_*$, and $0<\delta<1/2$.
There is
$\psi\in\PSH(X,\omega)\cap C^0(X)$, smooth off $E$, such that
$\omega_\beta=\omega+i\partial_X\bar\partial_X\psi$
is a K\"ahler cone metric of angle $2\pi\beta$. For this background,
there are nonnegative $u_0,u_1\in C_c^\infty(N)$ satisfying
\[
 \lambda\omega_\beta\le\omega_\beta+i\partial_X\bar\partial_Xu_j
 \le\Lambda\omega_\beta,\qquad j=0,1,
\]
on $X\setminus E$, for some $0<\lambda\le\Lambda$.

Let $\varphi_t$ be the weak geodesic relative to $\omega$ joining
$\psi+u_0$ to $\psi+u_1$. There are an open neighborhood $V\Subset N$
of $x_*$ and $v\in C^\infty(V)$ such that $\varphi_t=v$ near
$V\times[\delta,1-\delta]$. In holomorphic coordinates
$(z,w_1,\ldots,w_{n-1})$ with $E=\{z=0\}$, we have
$\psi=\varepsilon|z|^{2\beta}$ on $V$, with $\varepsilon>0$, and
\begin{equation}\label{eq:conic-main-rank}
 \omega+i\partial_X\bar\partial_X\varphi_t
 =\lambda_0 i\sum_{\alpha=1}^{n-1}
       dw_\alpha\wedge d\bar w_\alpha,
 \qquad \lambda_0>0,
\end{equation}
on $V\times[\delta,1-\delta]$. Thus the spatial form has rank $n-1$
on $V$, with null direction $\partial/\partial z$. When $n=1$, the
sum is empty and the spatial form is zero.
\end{maintheorem}

The constructed cone background and endpoint pair may depend on $\delta$.
Both endpoint metrics retain their cone angle. On the indicated region,
the intermediate potential $\varphi_t=v$ is smooth across $E$ and loses
its cone singularity; the relative potential $\varphi_t-\psi=v-\psi$
is not smooth across $E$. Remark~\ref{rem:conic-tangential} gives a
tangential version that retains the normal cone angle.

For the Poincar\'e case, let $\overline X$ be a compact K\"ahler
manifold of complex dimension $n\ge1$, let $E$ be a nonempty smooth
divisor, and set
$X=\overline X\setminus E$. In holomorphic coordinates
$(z,w_1,\ldots,w_{n-1})$ with $E=\{z=0\}$, $z$ is the normal
coordinate and $w$ gives the tangential coordinates. The local model is
\begin{equation}\label{eq:poincare-model}
 \frac{i\,dz\wedge d\bar z}{|z|^2(\log|z|^2)^2}
 +i\sum_{\alpha=1}^{n-1}dw_\alpha\wedge d\bar w_\alpha.
\end{equation}
It has potential $-\log(-\log|z|^2)+|w|^2$ for $0<|z|<1$.
The normal distance to $E$ is infinite, so a K\"ahler metric on $X$
uniformly comparable to this model near $E$ is complete.

Let $s_E$ be a defining section of $\mathcal O(E)$, equipped with
a smooth Hermitian metric $\mathfrak h$, and choose $\lambda$ large
enough that
\[
 s=\tfrac12(\lambda-\log|s_E|_{\mathfrak h}^2)>1
 \quad\text{on }X.
\]
Locally $s=-\log|z|+O(1)$, so $s\to\infty$ means approaching $E$.
For large $s$, the Poincar\'e distance to a fixed compact part of $X$
is comparable to $\log s$.

Fix a K\"ahler class on $\overline X$ and a point $x_*\in E$.
Balls $B_r$ are taken for a fixed smooth metric on $\overline X$.
For a Poincar\'e reference form $\omega_P$ in the sense of
\cite[Definition 0.1]{AuvrayPoincare}, let
$\mathcal H_{PM}(X,\omega_P)$ denote the bounded relative potential
class consisting of $u\in C^\infty(X,\R)$ such that
\[
 C^{-1}\omega_P\le\omega_P+i\partial_X\bar\partial_Xu
       \le C\omega_P,\qquad
 \sup_X|\nabla_{\omega_P}^ku|_{\omega_P}<\infty
       \quad(k\ge0),
\]
for some $C\ge1$. All derivatives and norms are taken with respect
to $\omega_P$. If $\rho_P$ is a local potential of $\omega_P$,
the corresponding total potential is $\rho_P+u$.

\begin{maintheorem}[Poincar\'e ends]\label{thm:poincare-boundary}
There exist a complete Poincar\'e reference form $\omega_P$ in the
restricted class and $u_0,u_1\in\mathcal H_{PM}(X,\omega_P)$ satisfying
\[
 |\nabla_{\omega_P}^ku_j|_{\omega_P}=O(s^{-N})
 \qquad(s\to\infty)
\]
for $j=0,1$, every integer $k\ge0$, and every $N>0$.
For every sufficiently small $r>0$, there are a nonempty open set
$V_r\subset X\cap B_r(x_*)$ and an open interval
$\tfrac12\in I_r\Subset(0,1)$ such that their weak geodesic $u_t$
relative to $\omega_P$ is independent of time on
$V_r\times I_r$. There, $\omega_P+i\partial_X\bar\partial_Xu_t$
has rank $n-1$, with kernel spanned by $\partial/\partial z$ in fixed
holomorphic coordinates $(z,w)$ near $x_*$ with $E=\{z=0\}$.
\end{maintheorem}

In the above theorem, we construct endpoints that decay faster than every power of $s^{-1}$ as $s\to\infty$. Even with this stronger endpoint decay, the
spatial form loses rank on open sets approaching $E$, including at
$t=1/2$. Thus rapid endpoint decay does not ensure second-order decay
along the weak geodesic. The intervals $I_r$ may shrink as $r\to0$.

\subsection{Mechanism and organization}

We construct endpoint potentials whose rooftop envelope is spatially
degenerate on a chosen region. Boundary contacts and harmonic measure
force the weak geodesic to agree with this envelope for a positive
time interval. In higher dimensions, we apply the same argument along
holomorphic curves.

Section~\ref{sec:preliminaries} gives the envelope, harmonic-measure
and gluing arguments. Section~\ref{sec:compact} gives the first
explicit construction in four parts: strictly subharmonic functions
on an annulus, model potentials on $\PP^1$, the construction on
K\"ahler products, and the extension to general compact K\"ahler
manifolds. Sections~\ref{sec:infinity} and
\ref{sec:divisor-singularities} treat the ALE and divisor examples,
respectively.
Appendix~\ref{sec:bounded-comparison} proves uniqueness of bounded
weak solutions and applies the comparison lemma in the ALE, conic,
and Poincar\'e settings.

\subsection*{Acknowledgements}

This work began during the second named author's first visit to IASM.
The authors thank Song Sun for his support.
The first named author thanks his supervisors, Julien Keller and Hugues Auvray,
for their guidance and encouragement.
This work was supported by the FRQNT grant
\emph{M\'etriques K\"ahl\'eriennes sp\'eciales singuli\`eres et non-compactes}
(\href{https://doi.org/10.69777/343263}{DOI: 10.69777/343263}).

\subsection*{Declaration of AI assistance}

Initial versions of the examples on $\PP^{1}$ and $\C$ were proposed by a
large language model. Building on these examples, the authors developed
a general framework for constructing weak K\"ahler geodesics with spatial
degeneracy and proved the results presented in this paper. AI also
assisted with language editing. The authors checked the mathematical
arguments, prepared the final exposition, and take full responsibility
for the mathematical content and presentation of the paper.

\section{Toolbox of construction}\label{sec:preliminaries}

We define the weak geodesic by envelope and establish the bounds used
in the constructions. We then express it through the shifted rooftop
envelopes and relate agreement with the rooftop to their one-sided
slopes. The envelope arguments apply to bounded continuous endpoints
on any K\"ahler manifold. We next prove the harmonic comparison
and gluing statements used in the examples.

\subsection{The weak geodesic by envelope}

We first specify the admissible functions and show that their envelope
solves the homogeneous complex Monge--Amp\`ere equation with the given
endpoints. We then compare this definition with the weak solutions
used in the geometric settings of the paper.

We write $\PSH(X,\omega)$ for the functions $f$ such that $\rho+f$
is plurisubharmonic in every chart with a smooth local potential $\rho$
of $\omega$. In particular these functions are upper semicontinuous and
$\omega+i\partial_X\bar\partial_X f\ge0$ as a current. 
For locally bounded potentials, the Monge--Amp\`ere measure is
understood in the sense of Bedford--Taylor
\cite{BedfordTaylor1976,BedfordTaylor1982};
see \cite[Section 2, Theorem 2.1]{BedfordTaylor1982} for the locally
bounded construction and its continuity under decreasing limits.

Let $u_0,u_1\in\PSH(X,\omega)\cap C^0(X)\cap L^\infty(X)$.
With $\Sigma$ and $\Omega$ as in the introduction, define
\[
\begin{split}
\mathcal B=\big\{W\in\PSH(X\times\Sigma,\Omega)\cap L^\infty:
\ \limsup_{(x',\tau')\to(x,\tau)}W(x',\tau')&\le u_j(x)\\
&\text{if }|\tau|=e^j,\quad j=0,1\big\}.
\end{split}
\]
The bounded Perron envelope is
\begin{equation}\label{eq:perron}
U=\left(\sup_{W\in\mathcal B}W\right)^*.
\end{equation}
All boundary limits in this definition are joint limits from inside
the product. The next lemma shows that $U$ is rotation invariant.
We call the resulting path
\[
 u_t(x)=U(x,\tau),\qquad t=\log|\tau|,
\]
the \emph{weak geodesic by envelope} joining $u_0$ to $u_1$.
Throughout the paper, we call this path simply the weak geodesic.
We recall the standard properties needed below.

\begin{lemma}\label{lem:perron}
The function $U$ is bounded and $\Omega$-plurisubharmonic, belongs to
$\mathcal B$, and attains the endpoint data. It is invariant under
rotation of $\tau$, and hence defines a path $u_t$.
If $A=\|u_1-u_0\|_\infty$, then
\begin{equation}\label{eq:barriers}
\max\{u_0-At,u_1-A(1-t)\}\le u_t\le(1-t)u_0+tu_1,
\end{equation}
and
\begin{equation}\label{eq:time-lipschitz}
|u_t(x)-u_s(x)|\le A|t-s|,\qquad s,t\in[0,1].
\end{equation}
Hence, $u_t(x)$ is Lipschitz in $t$ with a uniform Lipschitz bound $A$ for each fixed $x$.
Moreover, $(\Omega+i\partial\bar\partial U)^{n+1}=0$ on $X\times\Sigma$.
\end{lemma}

\begin{proof}
The maximum on the left of \eqref{eq:barriers} belongs to $\mathcal B$.
Restriction of any function in $\mathcal B$ to $\{x\}\times\Sigma$ is
subharmonic, so comparison with $(1-t)u_0(x)+tu_1(x)$ gives the upper
bound. This continuous upper bound persists under upper semicontinuous
regularization. The local upper-envelope theorem
\cite[Section 5]{BedfordTaylor1982}, applied after adding a local
potential of $\Omega$, gives $U\in\PSH(X\times\Sigma,\Omega)$.
The barriers give boundedness and joint endpoint attainment, hence
$U\in\mathcal B$.

Rotation of the annulus preserves $\mathcal B$, and therefore $U$.
Thus $u_t(x)=U(x,\tau)$ is well defined and convex in $t=\log|\tau|$.
For $0<s<t<1$, convexity gives
\[
\frac{u_s(x)-u_0(x)}s
\le\frac{u_t(x)-u_s(x)}{t-s}
\le\frac{u_1(x)-u_t(x)}{1-t}.
\]
The barriers bound the first quotient below by $-A$ and the last
above by $A$, proving \eqref{eq:time-lipschitz}, including the
endpoints by continuity in $t$.

Finally, choose a coordinate ball $B\Subset X\times\Sigma$ with a
smooth local potential $\rho$ of $\Omega$ on a neighborhood of
$\overline B$, and let $v=\rho+U$.
If $h\in\PSH(B)$ satisfies $h\le v$ outside a compact subset of $B$,
then $\max\{U,h-\rho\}$ on $B$, together with $U$ elsewhere,
defines a function in $\mathcal B$. The two formulas agree near
$\partial B$, and the endpoint limits are unchanged.
The envelope property therefore gives $h\le v$ on $B$.
Thus $v$ is maximal. The Bedford--Taylor characterization of maximal
locally bounded plurisubharmonic functions
\cite[Theorem 1]{Blocki2009Maximal} gives
$(i\partial\bar\partial v)^{n+1}=0$ on $B$, and hence
$(\Omega+i\partial\bar\partial U)^{n+1}=0$ on $X\times\Sigma$.
\end{proof}

\begin{remark}\label{rem:envelope-uniqueness}
On a general noncompact manifold, a bounded weak solution might be not
identified with the envelope by imposing the equation and
the two endpoint potentials; comparison at spatial infinity is needed.
In all the geometric examples in this paper, the bounded weak
solution is unique and agrees with~\eqref{eq:perron}.

For the compact, conic and Poincar\'e examples, existence follows from
\cite[Theorems 3--4]{Chen2000},
\cite[Theorem 2.1]{AleyasinChen2013}, and
\cite[Theorem 2.1]{AuvrayPoincare}, respectively. For ALE
spaces with $n\ge2$, we use \cite[Theorem D and Theorem 5.7]{YaoALE},
using the reference change \eqref{eq:ale-flat-reference} when needed.
The Euclidean case $n=1$ is treated in
Subsection~\ref{sec:noncompact}. These solutions have bounded potentials
relative to the reference forms used here.

For completeness, Lemma~\ref{lem:bounded-uniqueness} in
Appendix~\ref{sec:bounded-comparison} proves the identification
in all these cases using bounded comparison
\cite[Corollary 1.5]{KolodziejNguyen2023}.
\end{remark}

\subsection{Rooftop envelopes and their slopes}

We next prove the rooftop formula and identify exactly the times at
which the geodesic equals $P_0$. These facts reduce the later proofs
to estimates for $P_c$ as the real parameter $c$ varies.

For bounded continuous $\omega$-plurisubharmonic endpoints, define
\begin{equation}\label{eq:rooftop-intro}
P_c=P_\omega(u_0,u_1-c)
=\left(\sup\{f\in\PSH(X,\omega)\cap L^\infty(X):
 f\le u_0,\ f\le u_1-c\}\right)^*,\quad c\in\R.
\end{equation}
The star denotes upper semicontinuous regularization.
The envelopes $P_c$ in \eqref{eq:rooftop-intro} are bounded and
$\omega$-plurisubharmonic. A sufficiently negative constant is a
function in this class, and the continuous obstacles give a locally bounded
upper barrier. In particular, $P_c$ itself lies below both obstacles.

\begin{lemma}\label{lem:rooftop-properties}
For $c\le d$,
\begin{equation}\label{eq:rooftop-lipschitz}
P_c-(d-c)\le P_d\le P_c.
\end{equation}
For every fixed $x\in X$, the function $f_x(c)=P_c(x)$ is concave.
It has finite one-sided derivatives at every $c$, with
\[
-1\le f'_{x,+}(c)\le f'_{x,-}(c)\le0.
\]
\end{lemma}

\begin{proof}
Every function in the class defining $P_d$ also belongs to the class
defining $P_c$. Subtracting $d-c$ from a function in the class defining
$P_c$ gives one in the class defining $P_d$.
This proves \eqref{eq:rooftop-lipschitz}.
For $0\le a\le1$, the bounded $\omega$-plurisubharmonic function
$(1-a)P_c+aP_d$ lies below $u_0$ and below
$u_1-((1-a)c+ad)$. Hence
\[
(1-a)P_c+aP_d\le P_{(1-a)c+ad}.
\]
This is concavity. For each fixed $x$, concavity implies that the one-sided
difference quotients of $f_x$ are monotone, so their one-sided
limits exist. By \eqref{eq:rooftop-lipschitz}, these quotients
lie in $[-1,0]$, hence the limits are finite. Concavity also
implies that the right derivative is no larger than the left.
Therefore,
\[
-1\le f'_{x,+}(c)\le f'_{x,-}(c)\le0.
\]
\end{proof}

The compact rooftop formula is given in
\cite[Corollary 2.2]{DarvasRubinstein2016}. The following lemma proves
the formula for bounded envelopes on compact and noncompact manifolds,
using Kiselman's minimum principle
\cite{Kiselman1978,DarvasRubinstein2016}.

\begin{lemma}\label{lem:legendre}
The weak geodesic by envelope, $U(x, \tau) = u_t (x)$, satisfies
\[
P_c(x)=\inf_{0<t<1}\{u_t(x)-ct\},\qquad
u_t(x)=\sup_{c\in\R}\{P_c(x)+ct\}.
\]
\end{lemma}

\begin{proof}
Lift $U$ to the strip and apply Kiselman's minimum principle
\cite[Theorem 3.1 and Corollary 3.2]{DarvasRubinstein2016}.
In a local potential chart, the function
\[
\rho(x)+\inf_{0<t<1}(u_t(x)-ct)
\]
is plurisubharmonic. The infimum is globally bounded for each $c$.
Endpoint attainment gives
\[
\inf_{0<t<1}(u_t-ct)\le u_0,\qquad
\inf_{0<t<1}(u_t-ct)\le u_1-c.
\]
The infimum therefore belongs to the class defining $P_c$.
Conversely, we note that $F_c (x,\tau) = P_c (x) + c \log |\tau|$ 
is a function in $\mathcal B$. 
Hence, $F_c (x,\tau) \leq U(x,\tau)$, which gives $P_c(x) \leq u_t(x) - c t$. 
The first identity follows by taking the infimum over $t$.

For each fixed $x$ and $t\in(0,1)$, Jensen's inequality implies that the one-sided
difference quotients of $s\mapsto u_s(x)$ are monotone.
The one-side derivatives at each $t\in(0,1)$ exist and are bounded
by \eqref{eq:time-lipschitz}. Moreover, the left derivative is
no larger than the right derivative. Choosing $c$ between them,
we obtain
\[
u_s(x)-cs \ge u_t(x)-tc,\qquad 0\le s\le 1.
\]
For this choice of $c$, the first identity gives
\[
P_c(x)=\inf_{0<s<1}\{u_s(x)-cs\}=u_t(x)-ct.
\]
Since $P_c(x)+ct\le u_t(x)$ for every $c$, equality for this supporting
slope proves the second identity.

At the endpoints, recall $A=\|u_1-u_0\|_\infty$. The obstacle
inequalities give $P_{-A}=u_0$ and $P_A=u_1-A$, so the same formula
holds at $t=0,1$. Thus, the supremum is the pointwise value of $u_t$.
\end{proof}

\begin{lemma}\label{lem:exact-slopes}
For every $x\in X$,
\begin{equation}\label{eq:exact-time-interval}
\{t\in[0,1]:u_t(x)=P_0(x)\}
=[-f'_{x,-}(0),-f'_{x,+}(0)].
\end{equation}
\end{lemma}

\begin{proof}
Fix $x\in X$ and $t\in[0,1]$. By Lemma~\ref{lem:legendre},
the term $c=0$ gives $u_t(x)\ge P_0(x)$, and
\[
u_t(x)=P_0(x)
\quad\Longleftrightarrow\quad
P_c(x)+ct\le P_0(x)\quad\text{for every }c\in\R.
\]
For the finite concave function $f_x(c)=P_c(x)$, the last condition
says that the line of slope $-t$ through $(0,f_x(0))$ lies above
its graph. This holds if and only if
\[
f'_{x,+}(0)\le-t\le f'_{x,-}(0),
\]
or equivalently,
\[
-f'_{x,-}(0)\le t\le-f'_{x,+}(0).
\]
By Lemma~\ref{lem:rooftop-properties}, both one-sided derivatives
lie in $[-1,0]$, so this interval is contained in $[0,1]$.
This proves \eqref{eq:exact-time-interval}, including its endpoints.
\end{proof}
\subsection{Harmonic measure and exhaustion}\label{subsec:harmonic-measure}

We use harmonic measure to describe the time interval of spatial Hessian degeneracy. 
In this subsection, we introduce the concept of harmonic measure, and the basic facts.

\begin{definition}\label{def:harmonic-measure}
Let $B$ be a connected relatively compact domain in a Riemann surface,
with nonempty regular boundary. Thus every
continuous function $a$ on $\partial B$ has a unique harmonic extension
$H_a$ to $B$, continuous on $\overline B$ and equal to $a$ on $\partial B$.
For $z\in B$, the \emph{harmonic measure} $\omega_B^z$ is the probability
measure on $\partial B$ characterized by
\[
H_a(z)=\int_{\partial B}a(\xi)\,d\omega_B^z(\xi).
\]
\end{definition}

Here \emph{regular boundary} means that continuous boundary data are
attained at every boundary point by their harmonic extension. For
noncompact domains, this is a local condition at boundary points in the
ambient Riemann surface. All domains used for harmonic measure in the
constructions in this paper have regular boundary, including slit domains and
truncated domains.

In this subsection we restrict attention to continuous boundary functions
and step functions
\[
a=\sum_{j=1}^N\alpha_j\mathbf1_{E_j},
\]
where $\alpha_j\in\R$ and each $E_j$ is open or closed in the boundary.
Here $\mathbf1_E$ equals $1$ on $E$ and $0$ elsewhere.
For these data, let $H_a(z)=\int_{\partial B}a\,d\omega_B^z$.
For an open set $E\subset\partial B$, choose continuous functions
$0\le a_k\uparrow\mathbf1_E$. Their harmonic extensions increase and
satisfy $0\le H_{a_k}\le1$ by the maximum principle. Harnack's convergence
theorem says that an increasing sequence of harmonic functions bounded
above converges locally uniformly to a harmonic function
\cite[Theorem~3.8]{AxlerBourdonRamey2001}. Applying it in local coordinates,
we obtain the harmonic function
\[
\omega_B^z(E)=H_{\mathbf1_E}(z)=\lim_{k\to\infty}H_{a_k}(z),
\]
where the equality with harmonic measure follows by monotone convergence
in the defining integral. Closed sets are handled by complements.
Thus $H_a$ is harmonic for the boundary data considered here. Since $B$
has regular boundary, $H_a$ attains $a(\xi)$ at every continuity point
$\xi$ of $a$.

Let $D$ be a connected domain in a Riemann surface with regular boundary. If $D$ is
relatively compact, its harmonic measure is given by
Definition~\ref{def:harmonic-measure}. A bounded slit domain is allowed.
If $\overline D$ is not compact, we assume that $D$ has an increasing
exhaustion by domains $D_k$ with regular boundary and compact closure in
the ambient surface, such that
\begin{equation}\label{eq:exhaustion-assumption}
\partial D_k=A_k\cup\Gamma_k,\quad
A_k\subset\partial D,\quad \Gamma_k\subset D,\quad
A_k\uparrow\partial D.
\end{equation}
We take $A_k=\partial D_k\cap\partial D$, which is closed, and
$\Gamma_k=\partial D_k\cap D$, which is open in $\partial D_k$.
Thus $A_k$ is the part of the original boundary retained in
$\partial D_k$, and $\Gamma_k$ is the boundary introduced by truncation.
Neither is required to be a connected component of $\partial D_k$.
For a domain whose closure is not compact, we impose the additional
condition
\begin{equation}\label{eq:vanishing-harmonic-measure}
\omega_{D_k}^z(\Gamma_k)\longrightarrow0
\qquad\text{as }k\to\infty,\quad\text{for every fixed }z\in D.
\end{equation}
This is an additional restriction on the noncompact domain; it does
not follow from the existence of an exhaustion. For relatively compact
$D$ with regular boundary, take $D_k=D$ and $\Gamma_k=\varnothing$.

The next lemma defines harmonic measure on $D$ and gives the estimate
used below.

\begin{lemma}\label{lem:exhaustion}
Suppose $D$ has an exhaustion satisfying
\eqref{eq:exhaustion-assumption}. For an open or closed set
$E\subset\partial D$, the functions
\[
\omega_{D_k}^z(E\cap A_k)
\]
increase once $z\in D_k$, and their limit $\omega_D^z(E)$ is harmonic
in $z$. The limiting measure has total mass at most one.
If \eqref{eq:vanishing-harmonic-measure} also holds, its total mass
is one and the following comparison holds.
Let $a$ be a bounded continuous function or a step function on
$\partial D$, and let $r$ be subharmonic and bounded above on $D$.
If, for a fixed $z_0\in D$,
\[
\limsup_{D\ni w\to\xi}r(w)\le a(\xi)
\qquad\text{for }\omega_D^{z_0}\text{-almost every }\xi\in\partial D,
\]
then, for every $z\in D$,
\begin{equation}\label{eq:unbounded-comparison}
r(z)\le\int_{\partial D}a\,d\omega_D^z.
\end{equation}
\end{lemma}

\begin{proof}
First let $B$ be a relatively compact domain with regular boundary, and
let $r$ be subharmonic and bounded above on $B$. For $M>0$, the function
\[
b(\xi)=\max\left\{-M,\limsup_{B\ni w\to\xi}r(w)\right\}
\]
is bounded and upper semicontinuous on $\partial B$. Choose continuous
functions $b_j\downarrow b$. The maximum principle gives
$r\le H_{b_j}$, and bounded convergence gives $r(z)\le\int b\,d\omega_B^z$.
If the boundary upper limits of $r$ are at most $a$ almost everywhere
and $M\ge\|a\|_\infty$, then $b\le a$, $\omega_B^{z_0}$-almost everywhere. 
Let $E\subset\partial B$ be the exceptional set, with
$\omega_B^{z_0}(E)=0$ for some $z_0\in B$.
For every nonnegative continuous function $g$ on $\partial B$,
Harnack's inequality \cite[Theorem~3.6]{AxlerBourdonRamey2001} gives
\[
\int_{\partial B}g\,d\omega_B^z
\le C(z,z_0)\int_{\partial B}g\,d\omega_B^{z_0},
\]
with a constant independent of $g$. By the Riesz representation theorem,
this is an inequality of Borel measures:
\[
\omega_B^z\le C(z,z_0)\omega_B^{z_0}.
\]
In particular, $\omega_B^z(E)=0$.
Thus the boundary assumption holds $\omega_B^z$-almost everywhere
for every $z\in B$. For background on harmonic measure, see
\cite[Chapter 4]{Ransford1995} and
\cite[Chapter III]{GarnettMarshall2005}.
Since $M\ge\|a\|_\infty$, we have $b\le a$
$\omega_B^z$-almost everywhere, and hence
\begin{equation}\label{eq:hm-comparison}
r(z)\le \int_{\partial B}b\,d\omega_B^z
     \le \int_{\partial B}a\,d\omega_B^z.
\end{equation}
The same argument proves this comparison for every bounded Borel
boundary function $a$.

Now, we assume $D$ has an exhaustion $D_k$ of $D$ satisfying
\eqref{eq:exhaustion-assumption}.
To prove monotonicity of $\omega^z_{D_k} (E \cap A_k)$, fix $k<l$ and take a nonnegative continuous
function $a$ on $\partial D_l$. Its harmonic extension $H_a$ satisfies
\[
\int_{\partial D_l}a\,d\omega_{D_l}^z
=\int_{\partial D_k}H_a\,d\omega_{D_k}^z
\ge\int_{A_k}a\,d\omega_{D_k}^z,
\]
since $H_a\ge0$ and $H_a=a$ on $A_k$. As this holds for every
nonnegative continuous $a$, it is an inequality of measures. Hence
\[
\omega_{D_k}^z(E\cap A_k)
\le\omega_{D_l}^z(E\cap A_k)
\le\omega_{D_l}^z(E\cap A_l).
\]
Thus the restrictions of $\omega_{D_k}^z$ to $A_k$, extended by zero
to $\partial D$, increase as measures. Their limit defines
$\omega_D^z$, with total mass at most one. For fixed $E$, each term
is harmonic: $E\cap A_k$ is closed when $E$ is closed, and is a
difference of two closed sets when $E$ is open. Harnack convergence
gives harmonicity of the limit. To compare null sets at different points,
fix $z,z_0\in D$ and choose a finite Harnack chain with compact closure
in $D$ joining them. It is contained in $D_k$ for all sufficiently large
$k$, so the preceding measure domination holds with a constant
independent of $k$. Restricting to $A_k$ and passing to the increasing
limit gives
\[
\omega_D^z\le C(z,z_0)\omega_D^{z_0}.
\]
Thus the null sets of the limiting measures are independent of the pole.
Under \eqref{eq:vanishing-harmonic-measure},
\[
\omega_D^z(\partial D)
=\lim_{k\to\infty}\bigl(1-\omega_{D_k}^z(\Gamma_k)\bigr)=1.
\]

Finally, choose $M\ge\max\{\sup_Dr,\|a\|_\infty\}$. Apply the
bounded Borel comparison above on $D_k$, using $a$ on $A_k$
and $M$ on $\Gamma_k$. A set of zero limiting harmonic measure has
zero measure at every finite stage, by monotonicity. We obtain
\[
r(z)\le\int_{A_k}a\,d\omega_{D_k}^z
       +M\omega_{D_k}^z(\Gamma_k).
\]
The first term converges to $\int_{\partial D}a\,d\omega_D^z$, since
$a$ is bounded and the measures increase. The second tends to zero
by \eqref{eq:vanishing-harmonic-measure}, proving
\eqref{eq:unbounded-comparison}.
\end{proof}

In particular, a common contact set of positive harmonic measure
means that $\omega_D^z(C)>0$. It is a statement about the boundary
of the curve domain. It is not a statement about area in the ambient
manifold.

\begin{lemma}\label{lem:positive-arc}
Assume \eqref{eq:exhaustion-assumption} and
\eqref{eq:vanishing-harmonic-measure}. Let $C\subset\partial D$ be open
or closed and contain a nonempty smooth open arc $J$, relatively open
in $\partial D$. Then
\[
\omega_D^z(C)>0\quad(z\in D),\qquad
\lim_{D\ni z\to\xi}\omega_D^z(C)=1\quad(\xi\in J).
\]
For a slit, the limit holds from both sides.
\end{lemma}

\begin{proof}
Choose an interior point $\xi$ of the arc and a small one-sided domain
$B\subset D$ with regular boundary beside it. On $\partial B$, take
continuous data $0\le a\le1$, equal to one near $\xi$, supported
on the true-boundary arc, and zero on the artificial boundary.
Let $H_a$ be the harmonic extension on $B$. Extending $H_a$ by
zero to $D$ gives a bounded subharmonic function by the gluing
lemma. Its finite-boundary upper values are at most $\mathbf1_C$,
so harmonic comparison gives
\[
H_a\le\omega_D^\cdot(C)\le1.
\]
Since $H_a>0$ in $B$, the strong maximum principle gives
$\omega_D^\cdot(C)>0$ throughout $D$. Since $B$ has regular boundary,
$H_a(z)\to1$ as $z\to\xi$, proving the limit by squeezing.
For a slit, apply the same argument on each side.
\end{proof}
\subsection{Boundary contact and spatial degeneracy}\label{sec:contact}

In this subsection, we show that a time interval of rank loss can be obtained using harmonic measure on a
curve domain $D$. We call the part of $\partial D$ where both endpoints
agree with $v$ the \textit{common contact set}. 
Under the assumptions below, its harmonic measure gives the length of this guaranteed time interval when the ``rank loss'' happens.

Let $X$ be a Riemann surface with a K\"ahler form
$\omega$, and let
\[
u_0,u_1,v\in\PSH(X,\omega)\cap C^0(X)\cap L^\infty(X),
\qquad v\le u_0,u_1.
\]
Write $P_c=P_\omega(u_0,u_1-c)$ and let $u_t$ be the weak geodesic
by envelope connecting $u_0,u_1$.
Let $D\subset X$ be a connected domain, regular at its finite boundary
and satisfying
\eqref{eq:exhaustion-assumption}--\eqref{eq:vanishing-harmonic-measure}.
For the contact criterion, choose three disjoint boundary portions
$E_0,E_1,C$, each open or closed relative to $\partial D$, covering
$\partial D$ apart from finitely many junction points of zero harmonic
measure. These sets need not be connected or be connected components
of $\partial D$. Their roles are given by the contact equalities in
Theorem~\ref{thm:contact}.
When $\min(u_0,u_1)=v$ on $\partial D$, the exact contact partition is
\[
E_0=\{u_0=v<u_1\},\quad
E_1=\{u_1=v<u_0\},\quad C=\{u_0=u_1=v\},
\]
with all sets intersected with $\partial D$. The statements below
also allow additional contacts. Sections~\ref{sec:compact}--\ref{sec:divisor-singularities}
give the explicit choices of these sets in each construction.

The functions $m$ and $h$ are the harmonic measures of $E_1$ and
of the common contact set $C$, respectively:
\[
m(z)=\omega_D^z(E_1),\qquad h(z)=\omega_D^z(C),\qquad z\in D.
\]

\begin{theorem}
\label{thm:contact}
With the above notation and assumptions, suppose $v$ is smooth on $D$ and
\[
\omega+i\partial_X\bar\partial_Xv=0\quad\text{on }D.
\]
Assume that
\[
u_0=v\text{ on }E_0\cup C,
\qquad
u_1=v\text{ on }E_1\cup C.
\]
Suppose $h(z_*)>0$ at some point $z_*\in D$. Then, for every
$z\in D$ and $c\in\R$,
\begin{equation}\label{eq:main-rooftop}
P_c(z)-v(z)\le-cm(z)-\max(c,0)h(z),
\end{equation}
and
\begin{equation}\label{eq:main-identity}
u_t(z)=P_0(z)=v(z)
\qquad\text{if }m(z)\le t\le m(z)+h(z).
\end{equation}
On the nonempty open set
\[
\mathcal O=\{(z,t)\in D\times(0,1):m(z)<t<m(z)+h(z)\},
\]
the weak-geodesic potential agrees with $v$. In particular, it is
independent of time there and satisfies
\[
\omega+i\partial_X\bar\partial_Xu_t=0.
\]
The guaranteed interval of equality has length $h(z)>0$ at every
point $z\in D$.
\end{theorem}

\begin{proof}[Proof of Theorem~\ref{thm:contact}]
For each $c\in\R$, the bounded function $r_c=P_c-v$ is subharmonic
on $D$, since
\[
i\partial_X\bar\partial_Xr_c
=\omega+i\partial_X\bar\partial_XP_c\ge0.
\]
The inequalities $P_c\le u_0$ and $P_c\le u_1-c$, together with
continuity of $u_0,u_1,v$, give
\begin{equation}\label{eq:boundary-table}
\limsup_{D\ni z\to\xi}r_c(z)\le
\begin{cases}
0,&\xi\in E_0,\\
-c,&\xi\in E_1,\\
-\max(c,0),&\xi\in C.
\end{cases}
\end{equation}
Lemma~\ref{lem:exhaustion} gives \eqref{eq:main-rooftop}.

At $c=0$, this gives $P_0\le v$ on $D$. Since $v$ belongs to the class
defining $P_0$, the reverse inequality holds, so $P_0=v$ on $D$.
Since $r_c=P_c-v$ is subharmonic with boundary values comparison by \eqref{eq:boundary-table}, the maximum principle applies.
For $m(z)\le t\le m(z)+h(z)$,
\[
P_c(z)+ct-v(z)\le
\begin{cases}
c(t-m(z))\le0,&c\le0,\\
c(t-m(z)-h(z))\le0,&c\ge0.
\end{cases}
\]
Taking the supremum over $c$ and using the $c=0$ term proves
\eqref{eq:main-identity} by Lemma~\ref{lem:legendre}.

The harmonic functions $m,h$ are nonnegative and $m+h\le1$.
The strong maximum principle gives $h>0$ throughout $D$. Hence
$\mathcal O$ is nonempty and open by continuity. On $\mathcal O$,
$u_t=v$ is smooth and independent of time, and its spatial form is zero.
\end{proof}

\begin{remark}[Higher dimensional constructions]\label{rem:contact-higher-dimension}
Let $(X,\omega)$ be a K\"ahler manifold of dimension $n\ge2$, and let
\[
u_0,u_1,v\in\PSH(X,\omega)\cap C^0(X)\cap L^\infty(X),
\qquad v\le u_0,u_1.
\]
Let $D$ be a connected Riemann surface domain with the same regularity
and exhaustion assumptions as in Theorem~\ref{thm:contact}, and let
$B$ be a nonempty complex manifold of dimension $n-1$. Suppose an open
set $W\subset X$ is holomorphically identified with $D\times B$, and
this identification extends continuously to $\overline D\times B$.
Assume $v$ is smooth on $W$ and
\[
(\omega+i\partial_X\bar\partial_Xv)|_{D\times\{b\}}=0
\qquad(b\in B).
\]
For boundary portions $E_0,E_1,C$ as above, assume
\[
u_0=v\text{ on }(E_0\cup C)\times B,
\qquad
u_1=v\text{ on }(E_1\cup C)\times B.
\]
Here the boundary values are taken through the continuous extension
of the identification. Let $m(z)=\omega_D^z(E_1)$ and
$h(z)=\omega_D^z(C)$, and suppose $h(z_*)>0$ at one point of $D$.
Harmonic measure is taken on the curve domain $D$.

For the global rooftops $P_c=P_\omega(u_0,u_1-c)$ and the weak
geodesic by envelope $u_t$, one has, for every $(z,b)\in D\times B$,
\[
P_c(z,b)-v(z,b)\le-cm(z)-\max(c,0)h(z)\qquad(c\in\R),
\]
and
\[
u_t(z,b)=P_0(z,b)=v(z,b)
\qquad\text{if }m(z)\le t\le m(z)+h(z).
\]
On the nonempty open set
\[
\mathcal O=\{(z,b,t):z\in D,\ b\in B,\ m(z)<t<m(z)+h(z)\},
\]
the weak-geodesic potential agrees with $v$. In particular, it is
independent of time there and satisfies
\[
\omega+i\partial_X\bar\partial_Xu_t
=\omega+i\partial_X\bar\partial_Xv,
\qquad
\rank(\omega+i\partial_X\bar\partial_Xu_t)\le n-1.
\]
The guaranteed equality interval has length $h(z)>0$. No product
assumption on $\omega$ or independence of the potentials from $b$ is
required.

Indeed, fix $b\in B$. After adding a local potential $\rho$ of
$\omega$, the restriction of $\rho+P_c$ to the curve is subharmonic
and that of $\rho+v$ is harmonic. Thus $P_c(\cdot,b)-v(\cdot,b)$ is
bounded and subharmonic on $D$. The boundary comparison and rooftop
argument in the preceding proof apply to this function. Finally,
the smooth semipositive form of $v$ vanishes on the curve direction.
The Cauchy--Schwarz inequality makes its mixed terms with that direction
zero as well, so this direction lies in its kernel and the rank is
at most $n-1$.
\end{remark}

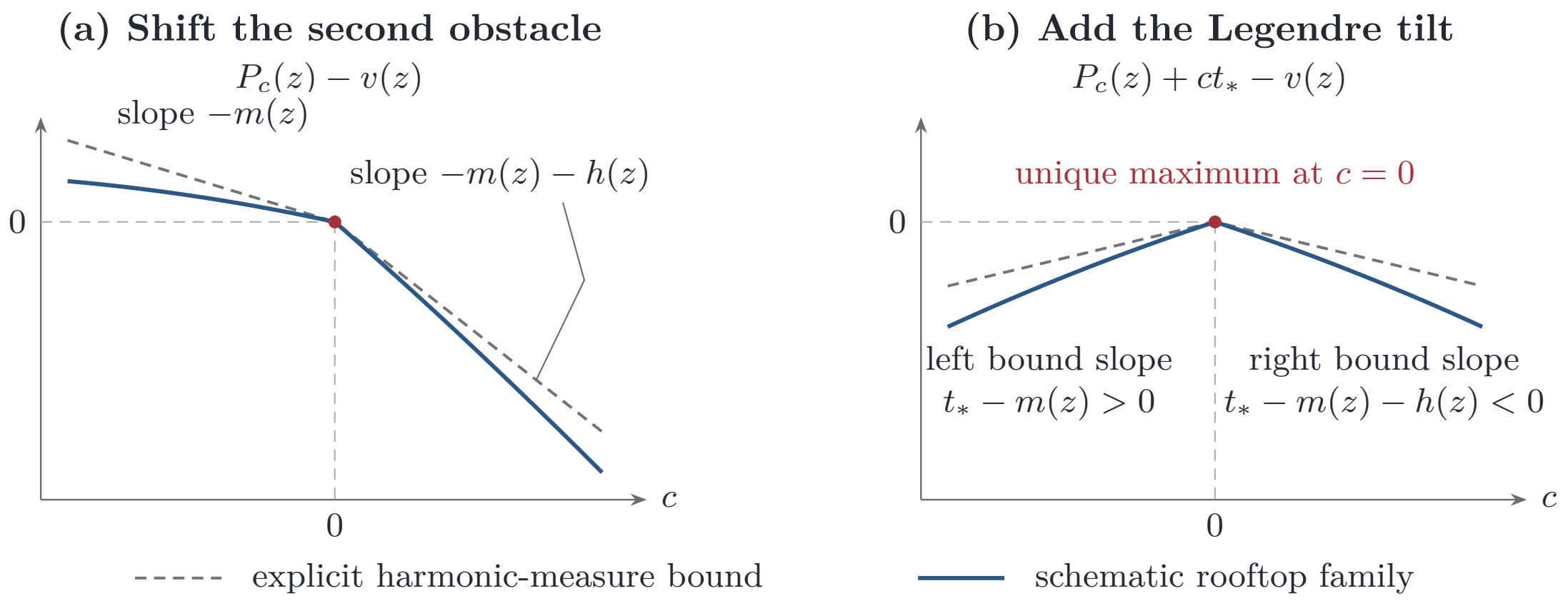
\begin{figure}[!htbp]
\centering
\begin{tikzpicture}[every node/.style={figtext}]
\path[use as bounding box] (-.2,-.75) rectangle (15.5,4.85);
\node[figtitle] at (3.55,4.55) {(a) Shift the second obstacle};
\node[figtitle] at (11.7,4.55) {(b) Add the Legendre tilt};
\node at (3.55,4.10) {\(P_c(z)-v(z)\)};
\node at (11.7,4.10) {\(P_c(z)+ct_*-v(z)\)};
\begin{scope}[shift={(3.6,2.78)},x=1.65cm,y=1.8cm]
 \draw[figaxis] (-1.65,-1.43)--(1.75,-1.43) node[right] {\(c\)};
 \draw[figaxis] (-1.65,-1.43)--(-1.65,.54);
 \draw[figguide] (-1.65,0)--(0,0);
 \node[left] at (-1.65,0) {\(0\)};
 \draw[figguide] (0,0)--(0,-1.43);
 \node[below] at (0,-1.43) {\(0\)};
 % The dashed majorant has illustrative m=.28 and h=.44.
 \draw[figbound] (-1.50,.42)--(0,0)--(1.50,-1.08);
 % Schematic rooftop: strictly below the majorant except at c=0.
 \draw[figcurve] plot[domain=-1.50:0,samples=40] (\x,{-.20*\x-.04*\x*\x});
 \draw[figcurve] plot[domain=0:1.50,samples=40] (\x,{-.80*\x-.04*\x*\x});
 \fill[figcontact] (0,0) circle[radius=1.7pt];
 \node[figlabel,anchor=south] at (-.68,.43) {slope \(-m(z)\)};
 \node[figlabel,anchor=south] at (.93,.12) {slope \(-m(z)-h(z)\)};
 \draw[figleader] (1.28,.10)--(1.40,-.30)--(1.13,-.814);
\end{scope}
\begin{scope}[shift={(11.75,2.78)},x=1.65cm,y=1.8cm]
 \draw[figaxis] (-1.65,-1.43)--(1.75,-1.43) node[right] {\(c\)};
 \draw[figaxis] (-1.65,-1.43)--(-1.65,.54);
 \draw[figguide] (-1.65,0)--(0,0);
 \node[left] at (-1.65,0) {\(0\)};
 \draw[figguide] (0,0)--(0,-1.43);
 \node[below] at (0,-1.43) {\(0\)};
 \draw[figbound] (-1.50,-.33)--(0,0)--(1.50,-.33);
 \draw[figcurve] plot[domain=-1.50:0,samples=40] (\x,{.30*\x-.04*\x*\x});
 \draw[figcurve] plot[domain=0:1.50,samples=40] (\x,{-.30*\x-.04*\x*\x});
 \fill[figcontact] (0,0) circle[radius=1.7pt];
 \node[figlabel,text=figcontact,anchor=south] at (0,.13) {unique maximum at \(c=0\)};
 \node[align=center] at (-.93,-.83) {left bound slope\\\(t_*-m(z)>0\)};
 \node[align=center] at (.95,-.83) {right bound slope\\\(t_*-m(z)-h(z)<0\)};
\end{scope}
\draw[figbound] (1.75,-.52)--(2.55,-.52);
\node[anchor=west] at (2.7,-.52) {explicit harmonic-measure bound};
\draw[figcurve] (9.0,-.52)--(9.8,-.52);
\node[anchor=west] at (9.95,-.52) {schematic rooftop family};
\end{tikzpicture}
\caption{Boundary comparison at \(z\in D\). Dashed curves show the bound
\eqref{eq:main-rooftop} and its shift by \(ct_*\); solid rooftop curves are
schematic. The slope labels refer to the dashed bounds.
For \(m(z)<t_*<m(z)+h(z)\), the maximum in (b) occurs at \(c=0\),
giving \(u_{t_*}(z)=P_0(z)=v(z)\).}
\label{fig:legendre-selection}
\end{figure}

Figure~\ref{fig:legendre-selection} separates the harmonic-measure
bounds from the actual rooftop graph.

\FloatBarrier

\begin{remark}\label{rem:main-slopes}
In this remark, we give a precise time interval of equality with the rooftop at a point.
For a fixed point $x$ (either $z$ or $(z,b)$), write $f_x(c)=P_c(x)$.
Concavity and the rooftop Lipschitz bound give finite one-sided
derivatives. Dividing the bound for $P_c-v$ by $c<0$ and $c>0$,
respectively, yields
\begin{equation}\label{eq:main-slopes}
f'_{x,-}(0)\ge-m(z),\qquad f'_{x,+}(0)\le-m(z)-h(z).
\end{equation}
Consequently $f'_{x,-}(0)-f'_{x,+}(0)\ge h(z)$, and
Lemma~\ref{lem:exact-slopes} gives the exact pointwise equality set
\[
\{t\in[0,1]:u_t(x)=P_0(x)\}
=[-f'_{x,-}(0),-f'_{x,+}(0)]\supset[m(z),m(z)+h(z)].
\]
The geometric interval is guaranteed, not necessarily maximal.
A strict corner characterizes positive-duration equality with $P_0$
at a point, not arbitrary spatial degeneracy. The open spatial
conclusion uses the continuous interval bounds and the smooth
comparison potential. No continuity of the actual rooftop derivatives
in $x$ is needed.
\end{remark}

\subsection{Gluing local potentials}\label{subsec:gluing}

In the constructions of Sections 3–5, we first define the data $(u_0,u_1,v)$ locally in a coordinate domain. 
The following gluing lemma allows us to extend these functions to $X$, preserving their values near a compact set 
and retaining global positivity.

Let $D\Subset X$ be a coordinate domain such that $\omega$ has a smooth
local potential on a neighborhood of $\overline D$. On that neighborhood,
let $u_0,u_1$ be smooth real functions and let $v$ be continuous and
$\omega$-plurisubharmonic, with
\[
 \omega+i\partial_X\bar\partial_Xu_j>0,\qquad
 v\le u_j,\qquad j=0,1.
\]
Fix a compact set $K\Subset D$ and a number $\eta>0$.

\begin{lemma}\label{lem:gluing}
Suppose
\begin{equation}\label{condi:gap-0}
 v>\eta\quad\text{on }K,\qquad
 u_j<-\eta\quad\text{near }\partial D.
\end{equation}
Then there are nonnegative functions $\widetilde u_j\in C_c^\infty(D)$ and a
bounded continuous $\omega$-plurisubharmonic function $\widetilde v$ on $X$,
all zero outside $D$, such that
\[
 \omega+i\partial_X\bar\partial_X\widetilde u_j>0,\qquad
 0\le\widetilde v\le\widetilde u_j.
\]
On a neighborhood of $K$ they satisfy
\[
 \widetilde u_j=u_j,\qquad \widetilde v=v.
\]
\end{lemma}

\begin{proof}
We will first define a precise regularized maximum of $\{u_j,0\}$ to obtain $\widetilde u_j$, and use
$\max\{v,0\}$ for $\widetilde v$. 
The strict gap from $0$, \eqref{condi:gap-0}, will give exact agreement near $K$ and vanishing near $\partial D$. 
The convexity and a derivative bound of the auxiliary function will preserve positivity of the endpoint K\"ahler forms.

Let $\kappa_\eta$ be an even smooth nonnegative mollifier of integral
one supported in $(-\eta,\eta)$. Define
\begin{equation}\label{eq:local-convex-cutoff}
 F_\eta(s)=\int_\R\max\{s-a,0\}\kappa_\eta(a)\,da.
\end{equation}
It's easy to see that $F_\eta$ is smooth, and its derivatives satisfy
\[
0 \le F'_{\eta}(s) \le 1, \qquad F''_{\eta} (s) \ge 0.
\]
$F_\eta$ is indeed a regularized maximum of $s$ and $0$; compare
\cite[Chapter I, Lemma 5.18]{DemaillyCADG}. It can be checked 
that $F_\eta(s)\ge\max\{s,0\}$, and the support condition gives
\[
 F_\eta(s)=0\text{ if }s\le-\eta,\qquad
 F_\eta(s)=s\text{ if }s\ge\eta.
\]

We now apply this function to the local potentials. On $D$, set
\begin{equation}\label{eq:local-global-potentials}
 \widetilde u_j=F_\eta(u_j),\qquad \widetilde v=\max\{v,0\}.
\end{equation}
Since $v\le u_j<-\eta$ near $\partial D$, all three functions vanish
near $\partial D$ and can extend outside $D$ by zero. 
The endpoint potentials are smooth with compact support in $D$.

To check positivity, choose a smooth local potential $\rho$ of $\omega$ near $\overline D$.
Write $\omega_{u_j}=\omega+i\partial_X\bar\partial_Xu_j$.
The chain rule gives
\begin{equation}\label{eq:local-cutoff-Levi}
\begin{split}
 \omega+i\partial_X\bar\partial_X\widetilde u_j
 ={}(1-F'_{\eta}(u_j))&\omega+F'_{\eta}(u_j)\,\omega_{u_j}\\
 &+F_\eta''(u_j)i\partial_Xu_j\wedge\bar\partial_Xu_j>0.
\end{split}
\end{equation}
The first line is a convex combination of positive definite forms;
the second is semipositive. Both $\rho+v$ and $\rho$ are
plurisubharmonic, so
\[
 \rho+\widetilde v=\max\{\rho+v,\rho\}
\]
is plurisubharmonic on $D$. Since $\widetilde v=0$ near $\partial D$,
extending $\widetilde v$ by zero gives a globally
$\omega$-plurisubharmonic function.
The global functions $\tilde u_j$ and $\tilde v$ are nonnegative by their definitions. 
Continuity and boundedness follow from the compact support and the definition of the functions in \eqref{eq:local-global-potentials}. 
The ordering follows from
\[
 0\le\widetilde v=\max\{v,0\}\le\max\{u_j,0\}
 \le F_\eta(u_j)=\widetilde u_j.
\]
Finally, continuity and $v>\eta$ on the compact set $K$ give
$u_j\ge v>\eta$ on a neighborhood of $K$. Then, $F_\eta(u_j)=u_j$ and
$\max\{v,0\}=v$, which complete the proof of the lemma.
\end{proof}
\section{The compact construction}\label{sec:compact}

We prove Theorem~\ref{thm:compact} in four parts. First, we construct
strictly subharmonic total endpoint potentials $q_0,q_1$ above a
harmonic function $p$ on an annulus, with prescribed contacts on
its boundary circles and a common slit. Second, we apply
Lemma~\ref{lem:gluing} to obtain model potentials on $\PP^1$.
Third, the contact theorem
gives spatial degeneracy on $\PP^1$ and its K\"ahler products.
Finally, we extend the planar functions to $\C$ and insert them into
a coordinate ball of an arbitrary compact K\"ahler manifold. We give the choices of
parameters and the positivity and contact checks explicitly.

\subsection{Strictly subharmonic functions on the annulus}
\label{subsec:compact-annulus}
\begin{figure}[!htbp]
\centering
\begin{tikzpicture}[every node/.style={figtext}]
\path[use as bounding box] (0,-.45) rectangle (15.6,5.65);
\node[figtitle] at (3.1,5.36) {(a) The slit annulus};
\node[figtitle] at (11.2,5.36) {(b) Logarithmic coordinates};
\begin{scope}[shift={(2.95,2.83)}]
 \fill[black!3,even odd rule] (0,0) circle[radius=1.78] (0,0) circle[radius=.72];
 \draw[figsecond,line width=1.05pt] (0,0) circle[radius=1.78];
 \draw[figfirst,line width=1.05pt] (0,0) circle[radius=.72];
 % Radii are compressed linearly in log radius, retaining the slit topology.
 \draw[figcontact,line width=1.8pt,line cap=round] (1.0733,0)--(1.4267,0);
 \fill[figcontact] (1.0733,0) circle[radius=1.35pt];
 \fill[figcontact] (1.4267,0) circle[radius=1.35pt];
 \node[align=center,text=figmuted] at (0,0) {outside\\\(A\)};
 \node at (-.45,1.24) {\(D=A\setminus\gamma\)};
 \draw[figleader] (-.51,-.51)--(-1.25,-1.86);
 \node[below,text=figfirst] at (-1.25,-1.88) {\(E_0:\ |z|=1\)};
 \draw[figleader] (1.25,-1.25)--(1.55,-1.86);
 \node[below,text=figsecond] at (1.55,-1.88) {\(E_1:\ |z|=e^3\)};
 \draw[figleader] (1.26,.03)--(1.73,.67)--(2.28,.67);
 \node[text=figcontact,anchor=west,figlabel] at (2.28,.67) {\(C=\gamma\)};
\end{scope}
\node[align=center,text=figmuted] at (6.70,2.62) {\(z=e^{s+i\theta}\)};
\draw[figaxis] (6.15,2.22)--(7.16,2.22);
\begin{scope}[shift={(8.10,1.27)},x=1.90cm,y=.48cm]
 \fill[black!3] (0,-0)--(3,0)--(3,6.2832)--(0,6.2832)--cycle;
 \draw[figguide] (0,0)--(3,0);
 \draw[figguide] (0,6.2832)--(3,6.2832);
 \draw[figfirst,line width=1.05pt] (0,0)--(0,6.2832);
 \draw[figsecond,line width=1.05pt] (3,0)--(3,6.2832);
 \draw[figcontact,line width=1.8pt,line cap=round] (1,3.1416)--(2,3.1416);
 \fill[figcontact] (1,3.1416) circle[radius=1.35pt];
 \fill[figcontact] (2,3.1416) circle[radius=1.35pt];
 \node[text=figcontact,above] at (1.5,3.45) {\(C:\ 1\le s\le2,\ \theta=0\)};
 \node at (1.5,1.55) {\(p=s\)};
 \node[left] at (0,6.2832) {\(\pi\)};
 \node[left] at (0,3.1416) {\(0\)};
 \node[left] at (0,0) {\(-\pi\)};
 \node[above] at (-.17,6.65) {\(\theta\)};
 \foreach \s in {0,1,2,3}{
   \draw[figguide] (\s,-.1)--(\s,-.28);
   \node[below] at (\s,-.31) {\(\s\)};
 }
 \node[anchor=west] at (3.14,-.51) {\(s\)};
 \draw[figmuted,-{Stealth[length=1.5mm]}] (1.24,6.2832)--(1.71,6.2832);
 \draw[figmuted,-{Stealth[length=1.5mm]}] (1.24,0)--(1.71,0);
 \node[above,text=figmuted] at (1.5,6.42) {top and bottom identified};
\end{scope}
\node[text=figfirst] at (2.6,-.15) {\(E_0:\ q_0=p<q_1\)};
\node[text=figsecond] at (7.75,-.15) {\(E_1:\ q_1=p<q_0\)};
\node[text=figcontact] at (12.6,-.15) {\(C:\ q_0=q_1=p\)};
\end{tikzpicture}
\caption{Contacts on the annulus: \(q_0=p\) on \(E_0\),
\(q_1=p\) on \(E_1\), and \(q_0=q_1=p\) on the slit \(C=\gamma\).
}
\label{fig:compact-contacts}
\end{figure}
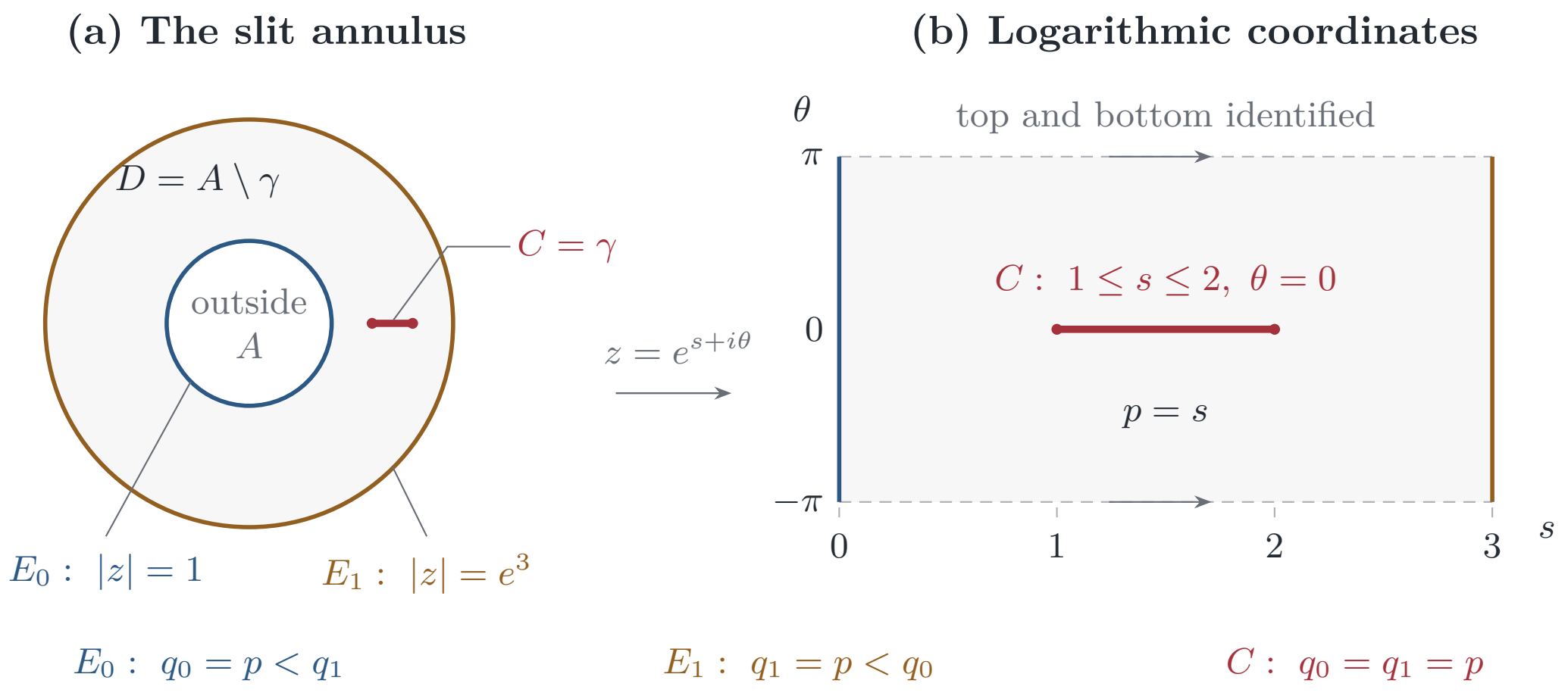

Throughout the section, we will use the coordinates $z=e^{s+i\theta}$, where $s=\log|z|$ is a spatial coordinate
and $t\in[0,1]$ denotes geodesic time.
Let
\[
 A=\{1<|z|<e^3\},\qquad p(s)=s,\qquad \gamma=[e,e^2].
\]
We seek smooth strictly subharmonic functions $q_0,q_1\ge p$ such
that $q_0=p$ on the inner circle, $q_1=p$ on the outer circle,
and both equal $p$ on the slit $\gamma$; see
Figure~\ref{fig:compact-contacts}. We obtain these contacts by adding
nonnegative terms to $p$.

Choose a smooth cutoff $\chi\in C^\infty(\R)$, with $0\le\chi\le1$,
such that
 $\chi(3-s)=\chi(s)$,
 $\chi=0$ outside $(1/2,5/2)$,and
 $\{\chi=1\}=[1,2]$.
We require
\[
 1-\chi(s)=e^{-1/(1-s)}
 \quad\text{for }s<1\text{ sufficiently close to }1,
\]
and use reflection near $s=2$. 

For $k>0$, let
\begin{equation}\label{eq:compact-B}
 B_0(s)=s^2e^{ks},\qquad B_1(s)=B_0(3-s).
\end{equation}
The factor $B_0$ vanishes at $s=0$ and $B_1$ at $s=3$, while
$1-\chi(s)\cos\theta$ vanishes exactly when $1\le s\le2$ and
$\theta=0$ modulo $2\pi$. Thus their products give the desired
circle and slit contacts. The next lemma shows that taking $k$
large also ensures strict subharmonicity.

\begin{lemma}\label{lem:compact-cutoff}
For all sufficiently large $k$, the functions in
\eqref{eq:compact-B} satisfy
\begin{equation}\label{eq:compact-laplacian-inequality}
 B_j''(s)>
 \left|(B_j\chi)''(s)-B_j(s)\chi(s)\right|,
 \qquad 0\leq s\leq3,\quad j=0,1.
\end{equation}
Consequently, for every $\varepsilon>0$, the functions
\begin{equation}\label{eq:compact-annular-endpoints}
 q_j(s,\theta)
 =s+\varepsilon B_j(s)\bigl(1-\chi(s)\cos\theta\bigr)
\end{equation}
are strictly subharmonic on $0\leq s\leq3$.
\end{lemma}

\begin{proof}
By symmetry, it suffices to consider $B_0$. The choice of $\chi$ gives
\begin{equation}\label{eq:compact-flat-cutoff-bound}
 |\chi'|+|\chi''|\leq C(1-\chi)^{3/4}.
\end{equation}
Indeed, near $s=1$ from the left, $1-\chi(s)=e^{-1/(1-s)}$,
so both derivatives of $\chi$ are bounded by
$2e^{-1/(1-s)}(1-s)^{-4}$. Note that
$e^{-1/(4(1-s))}(1-s)^{-4}$ is bounded as $s\to1^-$.
Reflection gives the same estimate near $s=2$. 
Let $C>0$ be an upper bound for
\[
 \frac{|\chi'(s)|+|\chi''(s)|}{(1-\chi(s))^{3/4}},
 \qquad 0\le s<1.
\]
The left-hand side of \eqref{eq:compact-flat-cutoff-bound} is zero outside $(0,1)\cup(2,3)$, so we prove the estimate.

On the support of $\chi$, we have $s\geq1/2$ and
\[
 \frac{B_0'}{B_0}=k+\frac2s\leq k+4,\qquad
 \frac{B_0''}{B_0}=k^2+\frac{4k}s+\frac2{s^2}\geq k^2.
\]
Then,
\begin{align*}
 \frac{B_0''-(B_0\chi)''+B_0\chi}{B_0}
 &=(1-\chi)\frac{B_0''}{B_0}
   -2\chi'\frac{B_0'}{B_0}-\chi''+\chi\\
 &\geq k^2(1-\chi)-C(k+1)(1-\chi)^{3/4}+\chi,\\
 \frac{B_0''+(B_0\chi)''-B_0\chi}{B_0}
 &=(1+\chi)\frac{B_0''}{B_0}
   +2\chi'\frac{B_0'}{B_0}+\chi''-\chi\\
 &\geq k^2-C(k+1).
\end{align*}
Young's inequality with exponents $4/3$ and $4$ gives
\[
 C(k+1)(1-\chi)^{3/4}
 \leq \frac{k^2}{2}(1-\chi)+C'k^{-2},\qquad k\geq1.
\]
For large $k$, the first lower bound is therefore at least
\[
 \frac{k^2}{2}(1-\chi)+\chi-C'k^{-2}
 \geq1-C'k^{-2}>0,
\]
and the second lower bound is also positive. Outside the support
of $\chi$, both numerators equal
$B_0''=e^{ks}(k^2s^2+4ks+2)>0$ for $0\le s\le3$, including at $s=0$.
This proves \eqref{eq:compact-laplacian-inequality}. Finally,
\[
 (\partial_s^2+\partial_\theta^2)q_j
 =\varepsilon\bigl[
 B_j''-\bigl((B_j\chi)''-B_j\chi\bigr)\cos\theta
 \bigr]>0.
\]
The Euclidean Laplacian in the coordinate $z=e^{s+i\theta}$ is
$\Delta_z=e^{-2s}(\partial_s^2+\partial_\theta^2)$, so this is
strict subharmonicity in the original complex coordinate.
\end{proof}

Fix $k$ as in Lemma~\ref{lem:compact-cutoff}. Since both $B_j$
are positive on $(0,3)$, \eqref{eq:compact-annular-endpoints} gives
\begin{equation}\label{eq:compact-contact-annulus}
 q_j\geq p,\qquad
 \{q_j=p\}\cap A
 =\gamma=\{e^s:1\leq s\leq2\}.
\end{equation}
At the boundary circles,
\[
 q_0=p<q_1\ \hbox{when }s=0,\qquad
 q_1=p<q_0\ \hbox{when }s=3.
\]
These contact properties hold for every $\varepsilon>0$. We now
choose $\varepsilon$ small enough to apply the gluing lemma.

\subsection{Model potentials on \texorpdfstring{$\PP^1$}{P1}}
\label{subsec:compact-model}

We first extend the annular formulas slightly past their boundary
circles. With $k$ fixed as above, choose $0<a\ll1$ and set
\[
 \widetilde A=\{e^{-a}<|z|<e^{3+a}\}.
\]
For $z=e^{s+i\theta}\in\widetilde A$, with $-a<s<3+a$, define
\[
\begin{aligned}
 p(z)&=s,\\
 q_0(z)&=s+\varepsilon s^2e^{ks}
              \bigl(1-\chi(s)\cos\theta\bigr),\\
 q_1(z)&=s+\varepsilon(3-s)^2e^{k(3-s)}
              \bigl(1-\chi(s)\cos\theta\bigr).
\end{aligned}
\]
These expressions are $2\pi$-periodic in $\theta$, so they define
smooth functions on $\widetilde A$. They agree with the functions
of Subsection~\ref{subsec:compact-annulus} on $A$, since the same
cutoff $\chi$ and the same formulas are used.

Since $\chi=0$ near $s=0,3$ and $B_j''>0$ at these points,
Lemma~\ref{lem:compact-cutoff} and continuity show that, for
sufficiently small $a>0$ independent of $\varepsilon$, the functions
$q_j$ are strictly subharmonic near $\overline{\widetilde A}$ for
every $\varepsilon>0$. The formulas also give $p\le q_j$ there.

Use a rescaled affine coordinate on $\PP^1$ in which
\[
 \rho_{\rm FS}(z)=\log(1+e^{-3}|z|^2)-c,\qquad
 \omegaFS=i\partial_X\bar\partial_X\rho_{\rm FS}.
\]
The constant $c$ does not change the reference form. To choose it,
write $p-\rho_{\rm FS}=f(s)+c$, where
\[
 f(s)=s-\log(1+e^{2s-3}).
\]
Then $f(3-s)=f(s)$ and
$f''(s)=-4e^{2s-3}/(1+e^{2s-3})^2<0$. In particular,
\[
 f(0)=f(3)>f(-a)=f(3+a).
\]
Choose $c=-(f(0)+ f(-a))/2$. Since $f(s)\ge f(0)$ on $[0,3]$,
we can choose $\eta>0$ sufficiently small, for instance, $\eta < (f(0)-f(-a))/8$, so that the local relative
potential $v^{\rm loc}=p-\rho_{\rm FS}=f(s)+c$ satisfies
\[
 v^{\rm loc}>4\eta\quad\text{on }\overline A,\qquad
 v^{\rm loc}<-4\eta\quad\text{on }\partial\widetilde A.
\]
Consequently, $v^{\rm loc}<-3\eta$ in a fixed neighborhood of
$\partial\widetilde A$. Since $B_j(s)(1-\chi(s)\cos\theta)$ is
bounded on $\overline{\widetilde A}$, we can choose $\varepsilon>0$
sufficiently small so that the local endpoint potentials
$u_j^{\rm loc}=q_j-\rho_{\rm FS}$ satisfy
\[
 v^{\rm loc}\le u_j^{\rm loc},\qquad
 u_j^{\rm loc}<-2\eta\quad\text{near }\partial\widetilde A.
\]
Their total forms are strictly positive, since
$\omegaFS+i\partial\bar\partial u_j^{\rm loc}
=i\partial\bar\partial q_j>0$, whereas
$\omegaFS+i\partial\bar\partial v^{\rm loc}=0$.

Lemma~\ref{lem:gluing}, with $K=\overline A$, now gives global
$u_0,u_1\in C^\infty(\PP^1)$ and a bounded continuous
$\omegaFS$-plurisubharmonic function $v$, all zero outside
$\widetilde A$, such that
\[
 \omegaFS+i\partial_X\bar\partial_Xu_j>0,\qquad 0\le v\le u_j.
\]
More explicitly, using $F_\eta$ from \eqref{eq:local-convex-cutoff},
these functions are
\[
 u_j=F_\eta(u_j^{\rm loc}),\qquad
 v=\max\{v^{\rm loc},0\}\quad\text{on }\widetilde A,
\]
and zero outside $\widetilde A$. The negative boundary gaps make
them zero on collars of both boundary circles, so the extensions
are well defined at $0$ and at $\infty$. The strict endpoint
positivity follows from the convex-combination formula
\eqref{eq:local-cutoff-Levi}; the maximum defining $v$ preserves
$\omegaFS$-plurisubharmonicity.
Near $\overline A$, we have $u_j=q_j-\rho_{\rm FS}$ and
$v=p-\rho_{\rm FS}$, so all circle and slit contacts are preserved.

\subsection{Construction on K\"ahler products}
\label{subsec:compact-sphere}

We first consider the higher-dimensional case on $\PP^1\times Y$,
where $(Y,\omega_Y)$ is a compact K\"ahler manifold.

\begin{proposition}\label{prop:compact-products}
The smooth endpoints $u_0,u_1$ on $\PP^1$ constructed above have the
following property. For every $0<\delta<1/2$, there is a nonempty
open set $V_\delta\Subset\PP^1$ such that their weak geodesic is
smooth and independent of time on a neighborhood of
$V_\delta\times[\delta,1-\delta]$, and
\[
 \omegaFS+i\partial_X\bar\partial_Xu_t=0
\]
there. The same endpoint pair works for every $\delta$.

For any compact K\"ahler manifold $(Y,\omega_Y)$, the pulled-back
endpoints on $\PP^1\times Y$, with reference form
$\pi_1^*\omegaFS+\pi_Y^*\omega_Y$, give a weak geodesic whose spatial
form is $\pi_Y^*\omega_Y$ on a neighborhood of
$V_\delta\times Y\times[\delta,1-\delta]$. Its rank there is
$\dim_\C Y$.
\end{proposition}

\begin{proof}
Use the global potentials $u_0,u_1,v$ constructed in the preceding
subsection. Let $D=A\setminus\gamma$, with inner circle $E_0$,
outer circle $E_1$, and common slit $C=\gamma$. The domain is
connected and regular for the Dirichlet problem, since every
connected component of its boundary contains more than one point;
see \cite[Chapter III, Corollary 6.2]{GarnettMarshall2005}.
Since $D$ is relatively compact,
the exhaustion conditions of Subsection~\ref{subsec:harmonic-measure}
hold with $D_k=D$ and no artificial boundary.

Let
\[
 m(z)=\omega_D^z(E_1),\qquad h(z)=\omega_D^z(C).
\]
The three boundary sets are disjoint and compact, so their indicator
functions are continuous on $\partial D$. Their harmonic extensions
attain these boundary values. The strong maximum principle therefore
gives
\[
 0<m<m+h<1\quad\text{on }D,
\]
and, as $z$ approaches an interior point of the slit from $D$,
\[
 m(z)\longrightarrow0,\qquad h(z)\longrightarrow1.
\]

On $D$, the function $v=p-\rho_{\rm FS}$ is smooth and
$\omegaFS+i\partial_X\bar\partial_Xv=0$. The boundary contacts
give $u_0=v$ on $E_0\cup C$ and $u_1=v$ on $E_1\cup C$.
Thus Theorem~\ref{thm:contact} gives
\begin{equation}\label{eq:compact-sphere-identity}
 u_t(z)=v(z)\qquad
 \bigl(z\in D,\ m(z)\le t\le m(z)+h(z)\bigr).
\end{equation}
Given $0<\delta<1/2$, choose a nonempty open set
$V_\delta\Subset D$ near an interior slit point so that
$m<\delta$ and $m+h>1-\delta$ on its closure.
These strict inequalities make \eqref{eq:compact-sphere-identity}
hold on an open neighborhood of $V_\delta\times[\delta,1-\delta]$.
There the geodesic is smooth, independent of time, and has zero
spatial form. The endpoint pair was fixed before $\delta$ was chosen.

On $\PP^1\times Y$, write
$\omega_{\rm prod}=\pi_1^*\omegaFS+\pi_Y^*\omega_Y$ and use the
pulled-back potentials $\pi_1^*u_j,\pi_1^*v$. The endpoint forms are
\[
 \omega_{\rm prod}+i\partial_X\bar\partial_X(\pi_1^*u_j)
 =\pi_1^*(\omegaFS+i\partial\bar\partial u_j)+\pi_Y^*\omega_Y>0.
\]
The same formula with $v$ gives a semipositive form, so $\pi_1^*v$
is a bounded continuous lower potential. All contacts persist on
$D\times Y$, and this lower form vanishes on each curve
$D\times\{y\}$.
Remark~\ref{rem:contact-higher-dimension} gives equality of the
weak geodesic with $\pi_1^*v$ on a neighborhood of
$V_\delta\times Y\times[\delta,1-\delta]$. Its spatial form there
is $\pi_Y^*\omega_Y$, of rank $\dim_\C Y$.
\end{proof}

\subsection{General compact K\"ahler manifolds}
\label{subsec:compact-general}

We now insert the planar construction into
a coordinate ball of an arbitrary compact K\"ahler manifold.
The reference form remains fixed, and the gluing lemma preserves
all contacts needed for the same harmonic-measure argument.

Using the potentials from Subsection~\ref{subsec:compact-model},
define on $\C$
\begin{equation}\label{eq:planar-seed}
 p=\rho_{\rm FS}+v,\qquad q_j=\rho_{\rm FS}+u_j,\qquad j=0,1.
\end{equation}
These extensions agree with the original functions near $\overline A$.
The function $p$ is continuous, subharmonic and bounded below;
the $q_j$ are smooth, strictly subharmonic and satisfy $p\le q_j$.

\begin{proof}[Proof of Theorem~\ref{thm:compact}]
Let $(X,\omega)$ and the nonempty open set $O\subset X$ be as in
the theorem. Choose holomorphic normal coordinates
$Z=(z,w)\in\C\times\C^{n-1}$ centered at a point of $O$, and a smooth
local potential $\rho$ for $\omega$ normalized so that
$\rho(Z)=|Z|^2+O(|Z|^3)$. Shrinking the
coordinate ball gives
\[
 \omega=i\partial_X\bar\partial_X\rho,\qquad
 \tfrac12|Z|^2\le\rho(Z)\le2|Z|^2
 \quad\text{on }\overline{\mathbb B_r}\Subset O.
\]
The potential is defined on a neighborhood of this closed ball.
When $n=1$, all transverse coordinates and factors below are omitted.

We add a common constant to the planar functions $p,q_j$ so that
$0\le p\le q_j$. This is possible because $p$ is bounded below, and
it changes neither their Hessians nor their contacts.
Choose $\ell>0$ and a transverse ball $B$ centered at zero so small
that
\[
 K=\ell\overline A\times\overline B\Subset\mathbb B_{r/4}.
\]
Only after fixing these spatial scales, define
\[
 M=\max\left\{1,\ \max_{j=0,1}
       \sup_{Z\in\overline{\mathbb B_r}}
            \bigl(q_j(z/\ell)+|w|^2\bigr)\right\},
 \qquad a=\frac{r^2}{16M}.
\]
The number $M$ is finite by continuity of $q_j$ and compactness of
$\overline{\mathbb B_r}$. On a neighborhood of $\overline{\mathbb B_r}$, put
\[
\begin{aligned}
 \widetilde q_j(Z)&=\frac{r^2}{4}
          +a\bigl(q_j(z/\ell)+|w|^2\bigr),\qquad j=0,1,\\
 \widetilde p(Z)&=\frac{r^2}{4}
          +a\bigl(p(z/\ell)+|w|^2\bigr).
\end{aligned}
\]
The endpoint total potentials are strictly plurisubharmonic:
writing $\zeta=z/\ell$, their Levi matrices are
\[
 \diag\left(
       \frac{a}{\ell^2}(q_j)_{\zeta\bar\zeta}(\zeta),
       a,\ldots,a\right)>0.
\]
The lower total potential $\widetilde p$ is continuous and
plurisubharmonic, and $\widetilde p\le\widetilde q_j$.

Set $\eta=r^2/16$. Since $K\Subset\mathbb B_{r/4}$, we have
$\sup_K\rho<r^2/8$, and hence
\[
 \widetilde p-\rho>\frac{r^2}{8}=2\eta
 \quad\text{on }K.
\]
On $\partial\mathbb B_r$, the choice of $a$ gives
\[
 \widetilde q_j-\rho
 \le\frac{r^2}{4}+\frac{r^2}{16}-\frac{r^2}{2}
 =-\frac{3r^2}{16}=-3\eta.
\]
By continuity, $\widetilde q_j-\rho<-\eta$ on a neighborhood of
$\partial\mathbb B_r$. Apply Lemma~\ref{lem:gluing} to the local
relative potentials $\widetilde q_j-\rho$ and
$\widetilde p-\rho$. It gives nonnegative endpoints
$u_j\in C_c^\infty(\mathbb B_r)\subset C_c^\infty(O)$ with
$\omega+i\partial_X\bar\partial_Xu_j>0$, and a global bounded
continuous $\omega$-plurisubharmonic function $0\le v\le u_j$.
On a neighborhood of $K$,
\[
 u_j=\widetilde q_j-\rho,\qquad
 v=\widetilde p-\rho,\qquad
 u_j-v=a(q_j-p)(z/\ell).
\]
Thus the inner-circle contact with $u_0$, the outer-circle contact
with $u_1$, and the common slit contact are all preserved throughout
the transverse ball $B$.

Let $D=A\setminus\gamma$, with harmonic measures $m,h$ as in the
preceding subsection. On $(\ell D)\times B$, the total potential
$\rho+v=\widetilde p$ is smooth and harmonic in $z$ for each fixed
$w$. Theorem~\ref{thm:contact} when $n=1$, and
Remark~\ref{rem:contact-higher-dimension} when $n>1$, therefore give
\[
 u_t(z,w)=v(z,w)
 \qquad\text{if }m(z/\ell)\le t\le m(z/\ell)+h(z/\ell).
\]
Fix a nonempty transverse ball $B'\Subset B$. For a given
$0<\delta<1/2$, choose $V_\delta\Subset D$ as in the preceding proof,
with $m<\delta$ and $m+h>1-\delta$ on its closure, and set
\[
 U_\delta=(\ell V_\delta)\times B'\Subset O.
\]
The strict harmonic-measure inequalities give $u_t=v$ on an open
neighborhood of $U_\delta\times[\delta,1-\delta]$. There the geodesic
is smooth and independent of time, and
\[
 \omega+i\partial_X\bar\partial_Xu_t
 =a\,i\sum_{\alpha=1}^{n-1}dw_\alpha\wedge d\bar w_\alpha.
\]
This form has rank exactly $n-1$; in dimension one the sum is empty
and the spatial form vanishes. All seed, scale and gluing parameters
were fixed before $\delta$ was chosen, so the same endpoint pair
works for every $\delta$.
\end{proof}
\section{Hessian degeneracy at infinity on ALE K\"ahler manifolds}
\label{sec:infinity}

We construct smooth endpoint potentials whose weak geodesic loses spatial
rank on open sets arbitrarily far out during interior time intervals.
We begin with rapidly decreasing endpoints on $\C$, extend the
construction to $\C^n$, and then obtain endpoints with polynomial decay
on the ALE manifolds in Theorem~\ref{thm:ale}. The main tools are the
boundary-contact criterion of Subsection~\ref{sec:contact}, which uses
harmonic measure to give a time interval of degeneracy, and smooth
cutoffs that preserve positivity. A change of reference form allows us
to carry the construction to an ALE end.

\subsection{The construction on \texorpdfstring{$\mathbb C^n$}{C\string^n}}
\label{sec:noncompact}

The construction has three steps. In the first step, we choose a narrowing
domain and construct two rapidly decreasing endpoint potentials with
boundary contacts determined by this domain. In the second step, we verify
the condition at infinity and apply the boundary-contact criterion to obtain
open regions and time intervals of spatial degeneracy. In the final step,
we extend the construction to higher dimensions using cutoffs and check
that positivity, decay, and rank loss are preserved.

\subsubsection{The planar domain and the endpoint potentials}\label{subsubsec:nc:planar}

\begin{figure}[!htbp]
\centering
\begin{tikzpicture}[every node/.style={figtext}]
\path[use as bounding box] (0,.35) rectangle (15.5,4.05);
% Schematic geometry; the text prescribes a smooth cap and exponential sides.
\def\hornpath{(3.1,3.25)
 .. controls (2.3,3.4248) and (.8,3.6) .. (.8,2.3)
 .. controls (.8,1.0) and (2.3,1.1752) .. (3.1,1.35)
 --plot[domain=3.1:13.5,samples=100] (\x,{2.3-.95*exp(-.23*(\x-3.1))})
 --plot[domain=13.5:3.1,samples=100] (\x,{2.3+.95*exp(-.23*(\x-3.1))}) --cycle}
\path[fill=figactive!8] \hornpath;
\draw[figguide] (.8,2.3)--(13.5,2.3);
\draw[figaxis] (13.5,2.3)--(15,2.3) node[right] {\(x\)};
\draw[figguide] (.8,1.5)--(.8,2.3);
\draw[figguide] (3.1,.75)--(3.1,3.25);
\node[below] at (.8,1.5) {\(R_0\)};
\node[below] at (3.1,.75) {\(R_0+2\)};
\draw[figcontact,line width=1.45pt] (3.1,3.25)
 .. controls (2.3,3.4248) and (.8,3.6) .. (.8,2.3)
 .. controls (.8,1.0) and (2.3,1.1752) .. (3.1,1.35);
\draw[figsecond,line width=1.1pt] plot[domain=3.1:13.5,samples=100]
 (\x,{2.3+.95*exp(-.23*(\x-3.1))});
\draw[figfirst,line width=1.1pt] plot[domain=3.1:13.5,samples=100]
 (\x,{2.3-.95*exp(-.23*(\x-3.1))});
\node[text=figcontact,anchor=west] at (.35,3.85) {\(\mathfrak C\): rounded cap};
\draw[figleader] (1.45,3.62)--(1.15,3.10);
\node[text=figsecond] at (6.3,3.65) {\(L_1:\ y=e^{-x}\)};
\draw[figleader] (6.3,3.36)--(5.9,2.80);
\node[text=figfirst] at (6.3,1.00) {\(L_0:\ y=-e^{-x}\)};
\draw[figleader] (6.3,1.27)--(5.9,1.80);
\node at (4.45,2.52) {\(D\)};
\node[text=figmuted] at (12.7,1.75) {towards infinity};
\end{tikzpicture}
\caption{The domain \(D\), with a rounded cap \(\mathfrak C\) on the left and
exponential sides \(L_0,L_1\) for \(x\geq R_0+2\).
The drawing is schematic.}
\label{fig:noncompact-horn}
\end{figure}
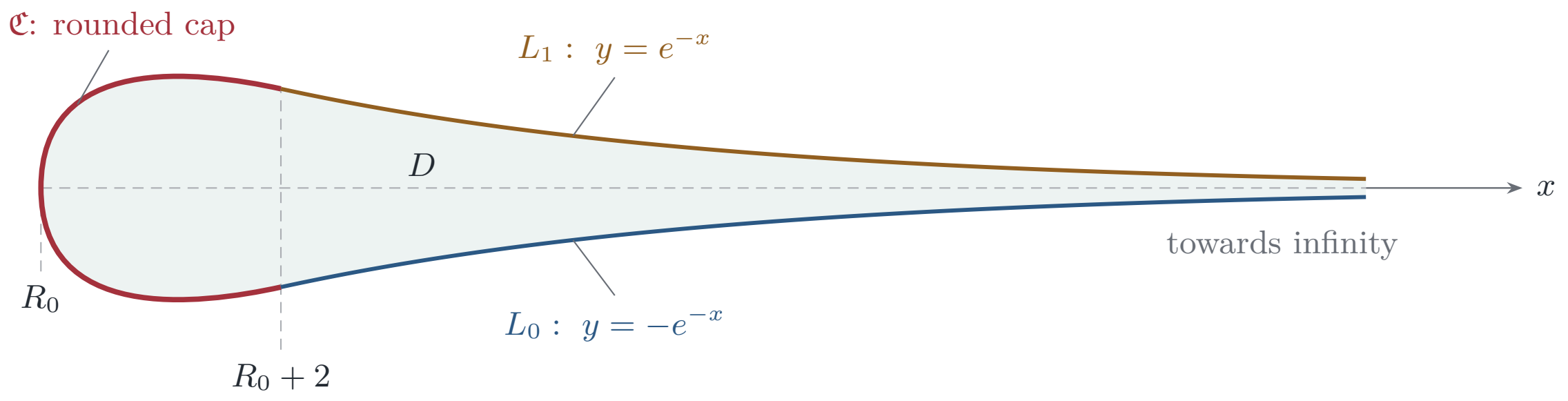

We work on \(\C\), write \(z=x+iy\) and
\(\Delta=\partial_x^2+\partial_y^2\), and set
\[
 \rho(z)=|z|^2,\qquad
 \omega_0=i\partial_X\bar\partial_X\rho.
\]
We would like to construct a weak geodesic whose total potential
\(\rho+u_t\) has zero spatial complex Hessian on open sets arbitrarily far out
in the narrow region between \(y=\pm e^{-x}\), during interior time
intervals. We first choose a smooth domain \(D\) which agrees with
this region for large \(x\); see Figure~\ref{fig:noncompact-horn}.

Fix a sufficiently large \(R_0>2\). For \(x\geq R_0+2\), let \(D\)
be the region between \(y=\pm e^{-x}\). Close its left end by a smooth
rounded cap in \(R_0\leq x\leq R_0+2\), joining the two curves smoothly
and symmetrically about the real axis. Denote the open lower and upper
sides by \(L_0,L_1\), and the cap, including the two junction points,
by \(\mathfrak C\). Thus \(D\) is simply connected, has smooth boundary, and lies
in a horizontal strip.

\FloatBarrier
Set $p(z)=\operatorname{Re}(z^2)=x^2-y^2.$ We will construct the total endpoint potentials by adding \(g_0,g_1\) to \(p\). 
Since \(p\) is harmonic, positivity of the endpoint forms reduces to positivity of \(\Delta g_j\).
To arrange the boundary contacts, we construct nonnegative functions \(g_0,g_1\): the first vanishes on the lower side and the cap, 
and the second on the upper side and the cap. 
Adding them to \(p\) will give the total endpoint potentials. 
We require positive Laplacians to obtain positive endpoint forms, and choose the behavior at infinity 
so that the endpoint potentials decay rapidly. 
We write \(\Sch(\C)\) for the smooth real functions whose derivatives of every order decay faster than every polynomial decay.

\begin{lemma}\label{lem:nc:zeros}
There are smooth functions \(g_0,g_1\geq0\) defined on $\C$, with $0<c \le \Delta g_j \le C$, whose boundary values satisfy
\begin{equation}\label{eq:nc:contacts}
 \begin{aligned}
 g_0&=0\quad\text{on }L_0\cup\mathfrak C,&\qquad g_0&>0\quad\text{on }L_1,\\
 g_1&=0\quad\text{on }L_1\cup\mathfrak C,&\qquad g_1&>0\quad\text{on }L_0.
 \end{aligned}
\end{equation}
Outside a compact set,
\[
 g_0=2(y+a)^2,\qquad g_1=2(y-a)^2,
\]
where \(a\in\Sch(\C)\) is nonnegative and satisfies
\(a(x,y)=e^{-x}\) for large \(x\) and \(|y|\leq2e^{-x}\).
\end{lemma}

\begin{proof}
On the unbounded end, the functions \(y+e^{-x}\) and \(y-e^{-x}\)
vanish on the lower and upper sides, respectively. We first replace
\(e^{-x}\) by a Schwartz function \(a\) on the plane, keeping its value
near these sides. Take \(\eta\in C^\infty(\mathbb R,[0,1])\), equal to zero for
\(x\leq R_0\) and one for \(x\geq R_0+1\), with derivatives bounded
independently of \(R_0\). Take
\(\vartheta\in C_c^\infty((-1,1),[0,1])\), equal to one on
\([-1/2,1/2]\). Define
\begin{equation} \label{def:afaroutpotential}
 a(x,y)=
 \begin{cases}
 e^{-x}\eta(x)\vartheta(y/x),&x>R_0,\\
 0,&x\leq R_0.
 \end{cases}
\end{equation}
The function is smooth across \(x=R_0\). It is supported where
\(x\geq R_0\) and \(|y|<x\); all its derivatives decay exponentially
there. It is therefore Schwartz. For large \(R_0\), \(|a_y|<1/2\).
The functions
\[
 f_{0,\infty}=y+a,\qquad f_{1,\infty}=y-a
\]
are consequently global submersions. Their zero curves equal the real
axis on the left and the respective graphs \(y=-e^{-x}\) and
\(y=e^{-x}\) on the right.

We next modify the two zero curves on a compact set so that each
contains the cap \(\mathfrak C\) and avoids the opposite side of \(D\). Complete
\(L_0\cup\mathfrak C\) to a smooth properly embedded line \(\Gamma_0\) by
continuing from the upper endpoint of \(\mathfrak C\) outside \(\overline D\),
above \(L_1\), then turning left above the whole bounded-height
strip and returning to the real axis sufficiently far to the left of
the cap, where \(a=0\). Complete
\(L_1\cup\mathfrak C\) similarly, with the unused branch below \(L_0\), to obtain
\(\Gamma_1\). Near either junction the unused branch can be separated
from the opposite side by a positive function flat at the endpoint.
Thus the continuation has the prescribed smooth jets, and
\[
 \Gamma_0\cap L_1=\varnothing,\qquad
 \Gamma_1\cap L_0=\varnothing.
\]
Each line can be chosen without self-intersections and to coincide with
\(f_{j,\infty}^{-1}(0)\) outside a large disc. Intersections between the
two unused branches cause no restriction.

To construct the defining functions, let \(Z_j=f_{j,\infty}^{-1}(0)\)
and first deform \(Z_j\) to \(\Gamma_j\), keeping both tails fixed.
Here only two compact simple arcs in a large disc need to be moved.
Keep smaller neighborhoods of their common end pieces fixed, and make
the remaining intersections transverse by a small perturbation. Deform
across the inner regions between the arcs. These
slides can be performed by smooth local graph deformations, so every
intermediate curve remains embedded.
Parametrize the resulting deformation by a smooth family of embeddings
\(\gamma_{j,s}:\R\to\R^2\), \(0\leq s\leq1\), with
$\operatorname{Im}\gamma_{j,0}=Z_j$, $\operatorname{Im}\gamma_{j,1}=\Gamma_j$, 
$\gamma_{j,s}(r)=\gamma_{j,0}(r)$ for $|r|\geq A$,
where \(A\) is fixed. Along each curve prescribe the velocity
\[
 V_{j,s}\bigl(\gamma_{j,s}(r)\bigr)
       =\partial_s\gamma_{j,s}(r).
\]
Using normal coordinates near the moving curves and a smooth cutoff,
extend these velocities to vector fields on \(\R^2\), smooth in
\((s,z)\) and supported in one fixed compact set. Their flows exist
for \(0\leq s\leq1\) and satisfy
\[
 \partial_s\Phi_{j,s}(z)=V_{j,s}\bigl(\Phi_{j,s}(z)\bigr),
 \qquad \Phi_{j,0}(z)=z.
\]
Uniqueness for this ODE gives
\(\Phi_{j,s}\circ\gamma_{j,0}=\gamma_{j,s}\), and hence
\(\Phi_{j,1}(Z_j)=\Gamma_j\). Thus
\[
 \Theta_j=\Phi_{j,1}^{-1},\qquad
 f_j=f_{j,\infty}\circ\Theta_j
\]
define smooth diffeomorphisms \(\Theta_j\), equal to the identity
outside a compact set, and smooth functions \(f_j\).
Then \(f_j^{-1}(0)=\Gamma_j\). It remains to obtain positive Laplacians
without changing these zero sets or the formulas at infinity. 

The square \(2f_j^2\) has the correct zero set, but
\[
 \Delta(2f_j^2)=4\bigl(|\nabla f_j|^2+f_j\Delta f_j\bigr)
\]
need not be positive. We first construct a nonnegative function
\(\varphi\), vanishing only at zero, so that
\(\varphi(f_j)\) has a uniformly positive Laplacian.
We then use a cutoff to restore \(2f_j^2\) at infinity.

All derivatives of \(f_j\) are bounded, and
\(|\nabla f_j|\geq a_0>0\). Outside a compact set, \(f_j=y\pm a\), so
\(\Delta(2f_j^2)=4+o(1)\) as \(|z|\to\infty\).
On the support of derivatives of \(a\), we have
\(|f_j|\leq x+O(1)\) and \(|\Delta f_j|=O(e^{-x})\). Hence, for
some \(C_2>0\),
\[
 \frac{|\Delta f_j|}{|\nabla f_j|^2}\leq \nu(|f_j|),
 \qquad \nu(s)=C_2e^{-s}\quad(s\geq0).
\]
The chain rule gives
\[
 \Delta\bigl(\varphi(f_j)\bigr)
 \geq|\nabla f_j|^2
 \bigl(\varphi''(f_j)
       -\nu(|f_j|)|\varphi'(f_j)|\bigr).
\]
For an even \(\varphi\), it therefore suffices to find $\varphi$ such that
\(\varphi''-\nu\varphi'\geq2\) and \(\varphi'>0\) on \((0,\infty)\).
Roughly, we set $w (s) = \varphi'(s)$. Solving $w' - \nu w = 4$ with $w(0)=0$, we can write the solution explicitly as:
\[
 w(s)=4\int_0^s e^{C_2(e^{-t}-e^{-s})}\,dt, \qquad \text{for } s\geq0.
\]
Note that $0 \le e^{-t} - e^{-s} \le 1$; then, $4s \leq w(s)\leq 4e^{C_2} s$. 
Since $(w-4s)' = \nu w$, we have
\[
w(s)- 4s = C_2 \int_0^s e^{-t} w(t)\,dt \le 4C_2e^{C_2} \int_0^s t e^{-t}\,dt.
\]
Hence, we have \(0\leq w(s)-4s\leq4C_2e^{C_2}\). 
To define a smooth function $\varphi$ in $\R$,
choose \(\varepsilon>0\) with \(s\nu(s)\leq1/2\) on
\([0,\varepsilon]\), and a smooth nondecreasing cutoff \(\theta\),
zero for \(s\leq\varepsilon/2\) and one for \(s\geq\varepsilon\).
Define
\[
 \varphi(s)=\int_0^{|s|}
       \bigl[4t+\theta(t)(w(t)-4t)\bigr]\,dt.
\]
This function is smooth and even, equals \(2s^2\) near zero, and is
positive away from zero. Note that
\[
 \varphi''-\nu\varphi'
 =4-4s\nu(s)(1-\theta)+\theta'(w-4s)\geq2.
\]
Thus, \(\Delta(\varphi(f_j))\geq2a_0^2\).
These Laplacians are also bounded above: \(\varphi''\) and
\(\varphi'(s)-4s\) are bounded, as are \(\nabla f_j\),
\(\Delta f_j\), and \(f_j\Delta f_j\).

Next, take a smooth cutoff
\(0\leq\chi_R\leq1\), equal to one on \(B_R\) and zero outside
\(B_{2R}\), with \(|\nabla\chi_R|=O(R^{-1})\) and
\(|\Delta\chi_R|=O(R^{-2})\), and define
\[
 g_j=(1-\chi_R)2f_j^2+\chi_R\varphi(f_j).
\]
On the annulus, $\{R<|z|<2R\}$,
\[
 \Delta g_j=(1-\chi_R)\Delta(2f_j^2)
       +\chi_R\Delta(\varphi(f_j))
       +2\nabla\chi_R\cdot\nabla (\varphi(f_j) - 2f_j^2)+(\varphi(f_j) - 2f_j^2)\Delta\chi_R.
\]
Note that $|(\varphi(f_j) - 2f_j^2) (z)|\leq C_0(1+|z|)$, $|\nabla (\varphi(f_j) - 2f_j^2) (z)|\leq C_0$; 
hence, the last two terms are \(O(R^{-1})\), whereas the sum of the first two
has a fixed positive lower bound for large \(R\).
Then, \(c_0\leq\Delta g_j\leq C_1\) on \(\C\) for some
\(c_0,C_1>0\). In conclusion, $g_j$, $j=0,1$ are the desired functions.
\end{proof}

We now define the total endpoint potentials \(q_j\) and the endpoint
potentials \(u_j\) relative to \(\omega_0\) by
\begin{equation}\label{eq:nc:endpoints}
 q_j=p+g_j,\qquad u_j=q_j-\rho=g_j-2y^2,\qquad j=0,1.
\end{equation}
Thus \(q_j\geq p\), with \(q_0=p\) on \(L_0\cup\mathfrak C\) and
\(q_1=p\) on \(L_1\cup\mathfrak C\). Since \(p\) is harmonic,
\[
 \omega_0+i\partial_X\bar\partial_Xu_j
 =i\partial_X\bar\partial_Xq_j
 =\tfrac14\Delta g_j\,\omega_0.
\]
The endpoint forms therefore have uniform positive lower and upper
bounds relative to \(\omega_0\). Outside a compact set, the construction
gives
\begin{equation}\label{eq:nc:tail-formulas}
 \begin{aligned}
 u_0&=4ya+2a^2,&\qquad u_1&=-4ya+2a^2.
 \end{aligned}
\end{equation}
Since \(a\) is Schwartz, these formulas show that
\(u_0,u_1\in\Sch(\C)\). This completes the endpoint construction.

\subsubsection{Hessian degeneracy on the narrowing domain}\label{subsubsec:nc:degeneracy}

We now show that the weak geodesic joining these endpoints has total
potential \(p\) on part of \(D\times(0,1)\).
Throughout this subsection, \(D\) denotes the domain constructed in 4.1.1 and shown in Figure 3.
We first construct a bounded function \(v\) whose total potential equals \(p\) near \(D\)
and which lies below both endpoints. We then verify the condition at
infinity and use the boundary contacts to determine the time intervals of equality.

To obtain a bounded relative potential, we replace the total potential
\(p\) by \(\max\{p,\rho-M\}\).
Choose \(M\) so large that
\[
 \rho-M\leq q_0,q_1\quad\hbox{on }\mathbb C,\qquad
 p= \rho-2y^2> \rho- M\quad\hbox{on }\overline D.
\]
The first condition follows from boundedness of the \(u_j\), and the
second from the strip containing \(D\). The function
\begin{equation}\label{eq:nc:v}
 v=\max\{p,\rho-M\}-\rho=\max\{-2y^2,-M\}
\end{equation}
is a bounded continuous \(\omega_0\)-plurisubharmonic function and
\(v\leq u_0,u_1\). It equals \(p-\rho\) on a neighborhood of
\(\overline D\). Thus \(\omega_0+i\partial_X\bar\partial_X v=0\) on \(D\), and
\eqref{eq:nc:contacts} gives exactly the boundary equalities required
in Theorem~\ref{thm:contact}.

Since \(D\) is unbounded, applying Theorem~\ref{thm:contact} also
requires \eqref{eq:vanishing-harmonic-measure}. We verify this by bounding
the harmonic measure of the artificial boundary of a truncation of
\(D\). For large \(T\), set
\[
 D_T=D\cap\{x<T\},\qquad
 \Gamma_T=\{T+iy:|y|<e^{-T}\}.
\]
These are bounded regular domains, with harmless corners at the
cross-section. Take \(W>\sup_{\overline D}|y|\) and let
\(\kappa=\pi/(2W)\). The positive harmonic function
\[
 B_T(x,y)=
 \frac{e^{\kappa(x-T)}\cos(\kappa y)}
      {\cos(\kappa e^{-T})}
\]
is at least one on \(\Gamma_T\) and nonnegative on the rest of
\(\partial D_T\). Hence
\begin{equation}\label{eq:nc:end-measure}
 0\leq\omega_{D_T}^{z}(\Gamma_T)\leq B_T(z)\longrightarrow0
 \qquad(T\to\infty)
\end{equation}
at every fixed \(z\in D\). Harmonic measure on \(\partial D\)
therefore has total mass one.

Each of \(L_0,L_1,\mathfrak C\) contains a nontrivial smooth boundary arc.
Lemma~\ref{lem:positive-arc} gives positive harmonic measure for each
part at every point of \(D\).
Define
\begin{equation}\label{eq:nc:mh}
 m(z)=\omega_D^z(L_1),\qquad h(z)=\omega_D^z(\mathfrak C).
\end{equation}
Then \(0<m(z)<m(z)+h(z)<1\). For \(P_c=P_{\omega_0}(u_0,u_1-c)\),
Theorem~\ref{thm:contact} therefore gives
\begin{equation}\label{eq:nc:envelope-bound}
 P_c(z)-v(z)\leq-cm(z)-\max(c,0)h(z),\qquad P_0=v
 \quad\hbox{on }D,
\end{equation}
and for \(z\in D\) and
\(m(z)\leq t\leq m(z)+h(z)\),
\begin{equation}\label{eq:nc:geodesic}
 u_t(z)=P_0(z)=v(z)=-2y^2,\qquad \rho(z)+u_t(z)=p(z).
\end{equation}
Indeed, \eqref{eq:nc:envelope-bound} bounds the left \(c\)-derivative
of \(P_c(z)\) at zero from below by \(-m(z)\), and the right derivative
from above by \(-m(z)-h(z)\). By concavity, \(c=0\) maximizes
\(P_c(z)+ct\) throughout this time interval.

Since \(D\) is unbounded and \(m,h\) are continuous with \(h>0\),
for every \(R>0\) there are a disc
\(B_R\Subset D\cap\{|z|>R\}\) and a nonempty open interval
\(I_R\Subset(0,1)\) such that \eqref{eq:nc:geodesic} holds on
\(B_R\times I_R\). The total spatial form is zero there.

\subsubsection{The construction in higher dimension}\label{subsubsec:nc:higher-dimension}

We keep the notation \(u_j\) for the planar endpoints on \(\mathbb C\)
and write \(\widetilde u_j\) for the endpoints constructed on
\(\mathbb C^n\), \(n\ge2\). The weak geodesic joining
\(\widetilde u_0,\widetilde u_1\) is denoted by \(\widetilde u_t\).

Pulling the planar endpoints back to \(\mathbb C^n\) does not give
Schwartz functions in all variables. We therefore use a cutoff in the
remaining variables. Write \(z=(z_1,z')\in\mathbb C\times\mathbb C^{n-1}\).
Let \(\chi\in C_c^\infty(\mathbb C^{n-1},[0,1])\) be one on the unit
ball and zero outside the ball of radius two. For a large fixed \(L\), set
\begin{equation}\label{eq:nc:lift}
 \chi_L(z')=\chi(z'/L),\qquad
 \widetilde u_j(z)=\chi_L(z')u_j(z_1),\qquad
 \widetilde q_j=|z|^2+\widetilde u_j.
\end{equation}
The functions \(\widetilde u_j\) are Schwartz because they have compact support in
\(z'\) and Schwartz decay in \(z_1\).

We check strict positivity before applying
Remark~\ref{rem:contact-higher-dimension}.
The planar positivity gives
\[
 (\widetilde q_j)_{z_1\bar z_1}
 =(1-\chi_L)+\chi_L(q_j)_{z_1\bar z_1}\ge c >0,
\]
with \(c\) independent of \(L\). All other entries of the complex
Hessian of \(\widetilde u_j\) are \(O(L^{-1})\), uniformly on
\(\mathbb C^n\), since \(u_j\) and its derivatives are bounded.
Thus, for sufficiently large fixed \(L\),
\[
 cI_n\le
 \bigl((\widetilde q_j)_{z_\alpha\bar z_\beta}\bigr)
 \le CI_n
\]
for some constants \(c,C>0\).

Let $\widetilde p(z)=\operatorname{Re}(z_1^2)+|z'|^2$. 
The cutoff preserves the inequality \(\widetilde q_j\ge \widetilde p\):
\begin{equation}\label{eq:nc:lift-lower}
 \widetilde q_j-\widetilde p
   =2(1-\chi_L)(y_1)^2+\chi_L(z')g_j (z_1)\geq0 .
\end{equation}
For sufficiently large \(M\), define
\[
 \widetilde v=\max\{\widetilde p,|z|^2-M\}-|z|^2
     =\max\{-2(y_1)^2,-M\}.
\]
As in \eqref{eq:nc:v}, \(\widetilde v\) is bounded, continuous, and
\(\omega_0\)-plurisubharmonic, with \(\widetilde v\le \widetilde u_0,\widetilde u_1\).
It therefore satisfies $ \widetilde v\le P_{\omega_0}(\widetilde u_0,\widetilde u_1) $.
On
\(D\times\{|z'|<L\}\) it is smooth and equals
\(-2( y_1)^2\). The spatial form is
\begin{equation}\label{eq:nc:rank}
 \omega_0+i\partial_X\bar\partial_X \widetilde v
       =i\sum_{\alpha=2}^n dz_\alpha\wedge d\bar z_\alpha .
\end{equation}
It vanishes along each curve \(D\times\{z'\}\). Since \(\chi_L=1\)
for \(|z'|<L\), the potentials \(\widetilde u_0,\widetilde u_1,\widetilde v\) are
independent of \(z'\) on \(D\times\{|z'|<L\}\). On \(\partial D\times B\), where
\(B=\{|z'|<L\}\), the boundary contact equalities follow from
\eqref{eq:nc:contacts} on each slice. Their harmonic measures are the
same functions \(m,h\) in \eqref{eq:nc:mh}, and
\eqref{eq:nc:end-measure} supplies the condition at infinity.
Remark~\ref{rem:contact-higher-dimension} therefore gives
\begin{equation}\label{eq:nc:lift-geodesic}
 \widetilde u_t(z)=-2( y_1)^2
 \quad\hbox{if }z_1\in D,\ |z'|<L,\
                m(z_1)\leq t\leq m(z_1)+h(z_1).
\end{equation}

\begin{proof}[The proof of the Euclidean case of Theorem~\ref{thm:ale}]
The construction is done in \ref{subsubsec:nc:planar} for $n=1$ and in \ref{subsubsec:nc:higher-dimension} for $n\ge2$. 
The region of rank loss is based on Theorem~\ref{thm:contact} and Remark~\ref{rem:contact-higher-dimension}, together with the discussion in \ref{subsubsec:nc:degeneracy} and \ref{subsubsec:nc:higher-dimension}. 
\end{proof}

The second named author asked whether the decay of the endpoint potentials and their
spatial derivatives at infinity is preserved by the weak geodesic at
every interior time. This construction gives a negative answer already
at second order: the endpoints are Schwartz, but the spatial complex
Hessian of the geodesic potential does not decay at infinity, even at
\(t=1/2\), as shown in Remark~\ref{rem:euclidean-midpoint} below.
For \(n\ge2\), the asymptotic behavior including all spatial derivatives,
for the corresponding HCMA problem with a disc parameter is proved in
\cite[Theorem B]{YaoALE}, under the ALE and boundary decay assumptions.

\begin{remark}\label{rem:euclidean-midpoint}
Complex conjugation preserves \(D\) and \(\mathfrak C\) and interchanges
\(L_0,L_1\). Hence, for every real \(x\in D\),
\[
\omega_D^x(L_0)=\omega_D^x(L_1),\qquad
m(x)=\frac{1-h(x)}2<\frac12<m(x)+h(x),
\]
since \(h(x)>0\). By continuity and
\eqref{eq:nc:geodesic}, \eqref{eq:nc:lift-geodesic}, the spatial
rank is \(n-1\) on open products near the real axis in \(D\), with
time intervals containing \(1/2\); in higher dimensions the spatial
sets are multiplied by \(B\). Such products occur arbitrarily far out.
\end{remark}
\subsection{The extension to ALE ends}
\label{sec:ale}

We extend the Euclidean construction under the holomorphic end and
decaying-potential hypotheses of Theorem~\ref{thm:ale}, written in
\eqref{eq:ale-potential}. The endpoint perturbations have polynomial
decay relative to the prescribed reference form. Relative to the flat
form on the end, their perturbations are rapidly decreasing. The proof
first makes a bounded change of reference form, then places the fixed
Euclidean horn in a cylinder embedded in the quotient end. A global
lower function makes the boundary-contact comparison applicable to
this same endpoint pair at every distance.

Recall that the end is biholomorphic to
$(\C^n\setminus\overline B_{R_0})/\Gamma$. In the holomorphic coordinates
$z$ on its covering space, with $r=|z|$, we assume
\[
 \omega=\omega_0+i\partial_X\bar\partial_X\psi,\qquad
 D^k\psi=O(r^{-\sigma-k})\quad(k\ge0),\qquad \sigma>0,
\]
where $\omega_0=i\partial_X\bar\partial_X|z|^2$ and $\psi$ is smooth,
real and $\Gamma$-invariant. Thus
$D^k(\omega-\omega_0)=O(r^{-\sigma-2-k})$. We construct endpoints
satisfying $D^ku_j=O(r^{-\sigma-k})$ for every $k\ge0$ and $j=0,1$;
all these derivatives are taken in the coordinates $z$.

\begin{proof}[Proof of Theorem~\ref{thm:ale}]
\emph{Step 1. Reduction to a Euclidean reference form.}
We first change the reference form so that it is Euclidean on the far
end. Increase $R_0$ if necessary, and choose a radial cutoff $\zeta_{R_0}$, zero for $r\leq R_0$ and one for
$r\geq2R_0$, with $D^k\zeta_{R_0}=O(R_0^{-k})$ in $\{R_0<r<2R_0\}$.
Extend $\Psi=\zeta_{R_0}\psi$  by zero over the compact part of $X$. 
This is a global smooth bounded function, since $\psi$ is invariant under $\Gamma$. 
Set
\begin{equation}\label{eq:ale-flat-reference}
 \widehat\omega=\omega-i\partial_X\bar\partial_X\Psi.
\end{equation}
It equals $\omega$ for $r\leq R_0$, $\omega_0$ for $r\geq2R_0$, 
and in $\{R_0<r<2R_0\}$, $\widehat\omega =\omega_0+O(R_0^{-\sigma-2})$.
For sufficiently large $R_0$, $\hat \omega$ is therefore K\"ahler, uniformly
equivalent to $\omega$. Since $\Psi$ is bounded,
subtracting $\Psi$ converts endpoint potentials and functions in the
bounded envelope class for $\widehat\omega$ into those for $\omega$. 
It therefore suffices to construct the endpoints and prove rank loss relative to $\widehat\omega$.

\emph{Step 2. Endpoint potentials on $\C^n/\Gamma$.}
Use the Euclidean endpoints $\widetilde u_j$ and the domain $D$ of
Subsection~\ref{sec:noncompact}. Write $z=(z_1,z')$, with $z_1=x+iy$.
Since $\Gamma$ acts freely on the sphere, choose a narrow cone
$\mathcal C_\Gamma$ about the positive real $z_1$-axis whose
$\Gamma$-translates are pairwise disjoint.

In the planar construction, choose
$\vartheta_\Gamma\in C_c^\infty((-1,1),[0,1])$, equal to one near
zero, with support so small that the corresponding cone in the
$z_1$-plane lies strictly inside $\mathcal C_\Gamma$. 
Define as \eqref{def:afaroutpotential} with the $\Gamma$ input:
\[
 a(x,y)=
 \begin{cases}
 e^{-x}\eta(x)\vartheta_\Gamma(y/x),&x>R_0,\\
 0,&x\leq R_0,
 \end{cases}
\]
where $\eta$ is the cutoff used in Lemma~\ref{lem:nc:zeros}.
All estimates remain valid after increasing the planar starting
radius. By \eqref{eq:nc:tail-formulas}, outside a compact set the
planar endpoint supports lie in this narrower cone. Choose $L$
sufficiently large for $D\subset\{|y|<L\}$,
and set
\[
 B=\{|z'|<L\},\qquad B'=\{|z'|<2L\}.
\]
The endpoint supports in the transverse directions lie in
$\overline{B'}$. In dimension one use the planar endpoints and omit
the transverse balls.

Translate $D$ and the endpoints far enough along the positive real
$z_1$-axis, keeping the same notation, so that the full endpoint
supports and $\overline D\times\overline{B'}$ lie inside
$\mathcal C_\Gamma$, with
\[
 D\subset\{x>R_1\},\qquad
 \widetilde u_j=0\quad\text{for }x\leq R_1,
\]
where $R_1$ is sufficiently large. Extend the endpoints to the
disjoint copies by
\begin{equation}\label{eq:ale-cylinder-endpoints}
 \widetilde u_j(\gamma z)=\widetilde u_j(z),
 \qquad z\in\mathcal C_\Gamma,\quad \gamma\in\Gamma,
\end{equation}
and set them to zero elsewhere. Since their supports lie inside
the cone, this gives smooth $\Gamma$-invariant potentials on
$\C^n/\Gamma$. Unitarity preserves the endpoint form bounds and
Schwartz decay. The functions vanish for $r\leq R_1$, since
$\Gamma$ preserves $r=|z|$. On $\mathcal C_\Gamma$, the original
inequality \eqref{eq:nc:lift-lower} and boundary contacts give
\begin{equation}\label{eq:ale-cutoff-inequality}
 2y^2+\widetilde u_j\geq0,
\end{equation}
with the prescribed equalities on $\partial D\times B$.
The construction allows arbitrarily large $R_1$, with uniform sup
norms and endpoint form bounds.

\emph{Step 3. Gluing to $X$ and proving rank loss.}
We first construct a bounded function $v$ below the endpoints,
equal to $-2y^2$ near $\overline D\times B$, as in
Subsection~\ref{subsubsec:nc:degeneracy}. The endpoint sup
norms and the height of $D$ are independent of $L$. Choose
$M>1+\max_j\|\widetilde u_j\|_\infty+2\sup_D y^2$ and define, on
$\mathcal C_\Gamma$,
\[
 v(z_1,z')=\max\{-2y^2,-M\}.
\]
By \eqref{eq:ale-cutoff-inequality} and the choice of $M$,
$v\leq\widetilde u_j$, and $v=-2y^2$ near
$\overline D\times B$, with exactly the boundary
contacts from Step 2. The total potential there is
$|z|^2+v=\operatorname{Re}(z_1^2)+|z'|^2$.

This local formula need not be $\Gamma$-invariant, so it does not
directly define a function on $(\C^n\setminus\overline B_{R_0})/\Gamma$.
In higher dimensions, $v$ also fails to equal $-M$ near the lateral
boundary of $\mathcal C_\Gamma$. We therefore apply a cutoff to $v+M$
and then subtract $M$. Use the same cutoff $\chi_L$ as in
\eqref{eq:nc:lift}, equal to one on $B$ and zero outside $B'$, and set
\[
 v_\chi=\chi_L(z')(v+M)-M
       =\chi_L(z')\max\{M-2y^2,0\}-M.
\]
We check positivity after this additional cutoff. In higher dimensions,
fix $M$, enlarge $L$ if necessary, choose $R_1$ as below, and then make
the final translation in Step~2. All derivatives below are Euclidean.
On $y^2<M/2$, the complex Hessian
$i\partial\bar\partial(|z|^2+v_\chi)$ is represented by the matrix
\[
 \begin{pmatrix}a&\mathbf b\\ \mathbf b^*&T\end{pmatrix},
\]
where
\[
 a=1-\chi_L,\qquad
 b_\beta=2iy(\chi_L)_{\bar z_\beta},\qquad
 T=I+(M-2y^2)(\chi_L)_{z'\bar z'}.
\]
The cutoff defined in \eqref{eq:nc:lift} satisfies
\[
 |\nabla_{z'}\chi_L|^2\le CL^{-2}(1-\chi_L),
\]
by Taylor's theorem applied to $1-\chi\ge0$ and scaling, with $C$
depending only on $\chi$. The same matrix represents the complex Hessian
of $|z|^2+\chi_L(M-2y^2)-M$ on $y^2<M$, where
\[
 |\mathbf b|^2=y^2|\nabla_{z'}\chi_L|^2\le CML^{-2}a,
 \qquad T\ge(1-CML^{-2})I.
\]
For sufficiently large $L$, $T\ge\tfrac12 I>0$ and
\[
 a-\mathbf bT^{-1}\mathbf b^*
 \ge a-2|\mathbf b|^2
 \ge a(1-C'ML^{-2})\ge0.
\]
Taking the maximum with $-M$ gives $v_\chi$ on this strip; since
$v_\chi=-M$ for $y^2\ge M/2$, it is $\omega_0$-plurisubharmonic
everywhere on the flat end, including across the maximum interface.
Moreover, $v_\chi\leq v\leq\widetilde u_j$, and
$v_\chi=v=-2y^2$ near $\overline D\times B$.
In dimension one take $\chi_L=1$ and omit the transverse variables.
Since $v_\chi=-M$ whenever $|y|\geq\sqrt{M/2}$ or $z'\notin B'$,
taking $R_1$ sufficiently large makes it equal to $-M$ near the
cone boundary for $r>R_1/2$. Copying it to the disjoint
$\Gamma$-translates and extending by $-M$ elsewhere gives a bounded
continuous $\omega_0$-plurisubharmonic function on this quotient end,
still denoted by $v_\chi$.

We now extend the data to $X$ by the regularized maximum construction in the
proof of Lemma~\ref{lem:gluing}. Fix $0<\eta<1$ and choose a smooth
radial function $b\in[-M,\eta]$, equal to $\eta$ on the compact part
and for $r\leq R_1/2$, and to $-M$ for $r\geq R_1$.
Its Hessian is $O(R_1^{-2})$, so taking $R_1>4R_0$
sufficiently large in Step 2 gives
\[
\widehat\omega+i\partial_X\bar\partial_Xb\geq\widehat\omega/2.
\]

Extend $\widetilde u_j$ by zero to $X$. With
$\operatorname{regmax}(s,t)=t+F_{\eta/2}(s-t)$ as in
\eqref{eq:local-convex-cutoff}, set
\[
 \widehat u_j=\operatorname{regmax}(\widetilde u_j,b),\qquad
 \widehat v=\operatorname{regmax}(v_\chi,b).
\]
Near the compact part, $b$ exceeds both $\widetilde u_j$ and
$v_\chi$ by at least $\eta$. Thus both regularized maxima
equal $b$ there and define global functions on $X$.
In particular, $\widehat v$ is bounded, continuous and
$\widehat\omega$-plurisubharmonic.
The endpoint form bounds follow from \eqref{eq:local-cutoff-Levi}.
For $r\geq R_1$, $\widehat u_j=\widetilde u_j$ because
$\widetilde u_j-b>1$. Near $\overline D\times B$ we also have
$v_\chi-b=M-2y^2>1>\eta/2$, so $\widehat v=v_\chi$ there.
Monotonicity of the regularized maximum therefore gives
\[
 \widehat v\leq\widehat u_j\quad\text{on }X,\qquad
 \widehat u_j=\widetilde u_j,\quad \widehat v=v_\chi=v=-2y^2
 \quad\text{near }\overline D\times B.
\]
All boundary contacts are preserved.

We can now apply Remark~\ref{rem:contact-higher-dimension}
(Theorem~\ref{thm:contact} when $n=1$).
The product $D\times B$ embeds holomorphically in the quotient end, the boundary contacts are
unchanged, and the spatial form of $\widehat v$ vanishes on each
curve $D\times\{z'\}$. The boundary regularity, positive common-contact
harmonic measure and exhaustion estimate \eqref{eq:nc:end-measure}
are preserved by translation. With $m,h$ denoting the harmonic
measures on this translated $D$, the weak geodesic satisfies
\[
 \widehat u_t(z_1,z')=\widehat v(z_1,z')=-2y^2
 \quad\text{for }m(z_1)\leq t\leq m(z_1)+h(z_1).
\]
On the strict-inequality region,
\begin{equation}\label{eq:ale-rank-identity}
 \widehat\omega+i\partial_X\bar\partial_X\widehat u_t
 =i\sum_{\alpha=2}^n dz_\alpha\wedge d\bar z_\alpha.
\end{equation}
Remark~\ref{rem:euclidean-midpoint} and continuity of $m,h$ give
open products containing time $1/2$ arbitrarily far out, for this
fixed endpoint pair.

Finally, return to the original reference form:
\begin{equation}\label{eq:ale-final-endpoints}
 u_j=\widehat u_j-\Psi,\qquad u_t=\widehat u_t-\Psi.
\end{equation}
Step 1 preserves the spatial forms, hence the rank identity.
Outside a compact set, $u_j=-\psi+\widetilde u_j$, which gives the
stated decay. The geodesic equals the smooth function
$\widehat v-\Psi$ on the identified open products.
\end{proof}

\section{Hessian degeneracy reaching singularities along divisors}
\label{sec:divisor-singularities}

We now consider endpoint metrics with cone or Poincar\'e singularities
along a smooth divisor. We prove Theorems~\ref{thm:conic-boundary}
and~\ref{thm:poincare-boundary} by constructing weak geodesics whose
spatial forms lose rank in the normal direction on open regions
approaching the divisor. In the cone case we use the compact
construction, adding the cone singularity to the endpoints while
preserving the boundary contacts. In the Poincar\'e case the metric
is complete on $\overline X\setminus E$, and we adapt the construction
at infinity from Section~\ref{sec:infinity}.

\subsection{Cone singularities}\label{sec:conic-boundary}

Let $(X,\omega)$ be a compact K\"ahler manifold of dimension $n\ge1$,
let $E\subset X$ be a nonempty smooth divisor, and fix $0<\beta<1$.
As in Theorem~\ref{thm:conic-boundary}, we write the cone background as
\[
 \omega_\beta=\omega+i\partial_X\bar\partial_X\psi,
 \qquad \psi\in\PSH(X,\omega)\cap C^0(X),
\]
where $\psi$ is smooth off $E$ and $\omega_\beta$ has cone angle
$2\pi\beta$ along $E$. Below we construct smooth endpoint potentials
$u_0,u_1$ so that $\psi+u_0,\psi+u_1$ define cone metrics, while their
weak geodesic relative to $\omega$ is smooth and has a null normal
direction on a neighborhood meeting $E$ during the prescribed interior
time interval. The endpoint cone angle is retained, but the intermediate
potentials lose their cone singularity on this neighborhood.

\begin{proof}[Proof of Theorem~\ref{thm:conic-boundary}]
\emph{Step 1. Smooth potentials and boundary contacts.}
Fix $x_*\in E$. We first construct $u_0,u_1,v$ using the slit-annulus
example in the normal coordinate near $x_*$. We place the origin inside
the slit domain so that its harmonic-measure interval contains
$[\delta,1-\delta]$, then use Lemma~\ref{lem:gluing} to extend the local
functions to $X$ while preserving their contacts.

Choose local holomorphic spatial coordinates $Z=(z,w)$ on $X$,
centered at $x_*$, with $E=\{z=0\}$ in this chart. Here $z\in\C$
is the normal coordinate and $w=(w_1,\ldots,w_{n-1})$ are tangential
coordinates. We normalize a local potential $\rho$ of $\omega$ while
preserving this equation for $E$. After subtracting its constant,
linear and pluriharmonic quadratic terms, its positive Hermitian
quadratic part can be made $|z|^2+|w|^2$. Shrinking the coordinate
ball gives
\[
 \omega=i\partial_X\bar\partial_X\rho,\qquad
 \tfrac12|Z|^2\le\rho(Z)\le2|Z|^2
 \quad\text{on }\overline{\mathbb B_r}\Subset N.
\]
When $n=1$, we omit the tangential coordinates and the transverse
factors below.

Recall from Section~\ref{sec:compact} the annulus
$A=\{\zeta\in\C:1<|\zeta|<e^3\}$ and the slit
$\gamma=[e,e^2]$ on the positive real axis.
On $D_*=A\setminus\gamma$, let $m,h$ be the harmonic measures of the
outer circle and slit, respectively. By
Subsection~\ref{subsec:compact-sphere}, $m\to0$ and $h\to1$ at an
interior slit point, so we can choose $\zeta_*\in D_*$ close to the slit with
\[
 m(\zeta_*)<\delta,\qquad m(\zeta_*)+h(\zeta_*)>1-\delta.
\]

We construct local total potentials by rescaling the planar data
$p,q_0,q_1$ in \eqref{eq:planar-seed}. Subtract a common constant so
that $0\le p\le q_0,q_1$, which is possible since $p$ is bounded below.
This preserves all contacts and Hessians. The function $p$ is harmonic
on $A$, and $q_0,q_1$ are smooth strictly subharmonic functions.
Our aim is to satisfy the hypotheses of Lemma~\ref{lem:gluing}: the
total potentials should lie above $\rho$ near the boundary contacts
and below $\rho$ near the coordinate-ball boundary. We apply the same
rescaling and add the same modifying terms to all three functions.

Place the scaled slit annulus in the $z$-coordinate by setting
\[
 z=\ell(\zeta-\zeta_*),\qquad
 D_\ell=\ell(D_*-\zeta_*),\qquad
 \ell(e^3+|\zeta_*|)<r/16.
\]
Thus $0\in D_\ell$, and the scaled annulus lies in $\{|z|<r/16\}$.
Choose $B=\{|w|<r_\perp\}$ with $r_\perp<r/16$, and set
\[
 K=\overline{D_\ell}\times\overline B\Subset\mathbb B_{r/2}.
\]
This compact set contains all the boundary contacts, and
$\sup_K\rho<r^2/64$.

With the radius $r$ fixed, choose
\[
 A_0=r^2/32,\qquad \lambda_0=1/32,\qquad a=\frac{A_0}{2M},
\]
where $M=\max\left\{1,\ \max_{j=0,1}
       \sup_{|\zeta|\le r/\ell+|\zeta_*|}q_j(\zeta)\right\}$.
On $\mathbb B_r$, with $\zeta=\zeta_*+z/\ell$, define
\[
\begin{aligned}
 \widetilde p&=A_0+a p(\zeta)+\lambda_0|w|^2,\\
 \widetilde q_j&=A_0+a q_j(\zeta)+\lambda_0|w|^2,
 \qquad j=0,1.
\end{aligned}
\]
These functions are respectively plurisubharmonic and strictly
plurisubharmonic, and $\widetilde p\le\widetilde q_j$.
Since $p\ge0$, we have $\widetilde p\ge A_0$, while
$\widetilde q_j\le A_0+A_0/2+\lambda_0r^2=5A_0/2$.
Consequently,
\[
 \widetilde p-\rho>A_0/2\text{ on }K,\qquad
 \rho-\widetilde q_j\ge 3A_0/2
       \text{ on }\{r/2\le|Z|\le r\}.
\]
Lemma~\ref{lem:gluing}, with $\eta=A_0/4$, gives nonnegative
$u_j\in C_c^\infty(N)$ with $\omega+i\partial_X\bar\partial_Xu_j>0$
and a bounded continuous $\omega$-plurisubharmonic function
$0\le v\le u_j$, all zero outside $\mathbb B_{r/2}$. Near $K$,
\[
 v=\widetilde p-\rho,\qquad u_j=\widetilde q_j-\rho,
 \qquad u_j-v=a(q_j-p)(\zeta).
\]
Thus all boundary contacts are preserved. Since $p$ is harmonic on
$A$, the function $v$ is smooth on $D_\ell\times B$, and its spatial
form vanishes on each normal slice with $w$ fixed. Moreover,
$\partial D_\ell\times\overline B$ is disjoint from $E$, because
$0$ lies in the interior of $D_\ell$. We use this separation in the
next step.

\emph{Step 2. Adding the cone singularity.}
We add a nonnegative cone potential $\psi$ to the endpoints, choosing
it to vanish near their boundary contacts. The lower function $v$
will remain unchanged.

To define the cone potential globally along $E$, let $s_E$ be the
defining section of $\mathcal O(E)$, which vanishes precisely on $E$.
Choose a smooth Hermitian metric $\mathfrak h$ with
$|e|_{\mathfrak h}^2=1$ near $\overline{\mathbb B_r}$ in the local
frame $s_E=ze$. Then $|s_E|_{\mathfrak h}^{2\beta}$ is globally
nonnegative, bounded and continuous, smooth off $E$, and equals
$|z|^{2\beta}$ near our construction.

Let $d=\operatorname{dist}(0,\partial D_\ell)>0$. Since the boundary
contacts are separated from $E$, choose a smooth $0\le\chi\le1$
equal to zero near $\partial D_\ell\times\overline B$, equal to one
on $\{|z|<d/3\}\cap\mathbb B_r$, and with
$\operatorname{supp}(1-\chi)\Subset\mathbb B_r\setminus E$.
In particular, $\chi=1$ near all of $E$. For $\varepsilon>0$ to be
chosen below, let
\[
 \psi=\varepsilon\chi|s_E|_{\mathfrak h}^{2\beta},\qquad
 \omega_\beta=\omega+i\partial_X\bar\partial_X\psi.
\]

Off $E$, direct differentiation gives
\[
\begin{aligned}
 i\partial_X\bar\partial_X\psi
 ={}&\varepsilon\beta^2|s_E|_{\mathfrak h}^{2\beta}
 i\partial_X\log|s_E|_{\mathfrak h}^2
 \wedge\bar\partial_X\log|s_E|_{\mathfrak h}^2\\
 &+\varepsilon\beta|s_E|_{\mathfrak h}^{2\beta}
 i\partial_X\bar\partial_X\log|s_E|_{\mathfrak h}^2\\
 &+\varepsilon i\partial_X\bar\partial_X
 \bigl((\chi-1)|s_E|_{\mathfrak h}^{2\beta}\bigr).
\end{aligned}
\]
The first term is nonnegative. The second has $\omega$-norm bounded
by a constant times $\varepsilon$, since the curvature of
$\mathfrak h$ is smooth and $|s_E|_{\mathfrak h}^{2\beta}$ is bounded.
The cutoff term has the same bound:
$(\chi-1)|s_E|_{\mathfrak h}^{2\beta}$ is globally smooth because
$\operatorname{supp}(\chi-1)$ is compact and disjoint from $E$.
Thus the sum of the last two terms lies between
$-\varepsilon C\omega$ and $\varepsilon C\omega$, with $C$
independent of $\varepsilon$. Hence
\[
 i\partial_X\bar\partial_X\psi\ge-\varepsilon C\omega.
\]
By continuity of $\psi$, this inequality extends across $E$ in the
sense of currents by the removable-singularity theorem
\cite[Chapter I, Theorem 5.24]{DemaillyCADG}.

By compactness, choose $0<\kappa\le1\le C_1$ such that
\[
 \kappa\omega\le\omega+i\partial_X\bar\partial_Xu_j\le C_1\omega.
\]
Fix $\varepsilon>0$ with $\varepsilon C<\kappa/2$. Then, on
$X\setminus E$,
\[
 \omega_\beta\ge\tfrac12\omega,\qquad
 \omega_\beta+i\partial_X\bar\partial_Xu_j\ge\tfrac{\kappa}{2}\omega.
\]
Using
$(\kappa-1)\omega\le i\partial_X\bar\partial_Xu_j\le(C_1-1)\omega$
also gives the required uniform comparison:
\[
 \tfrac{\kappa}{2}\omega_\beta
 \le\omega_\beta+i\partial_X\bar\partial_Xu_j
 \le2C_1\omega_\beta.
\]

Finally, $\psi\ge0$ gives $v\le u_j+\psi$, while $\psi=0$ near
$\partial D_\ell\times\overline B$. Thus the boundary contacts from
Step 1 are unchanged, and the lower function $v$ remains the same.

\emph{Step 3. Applying the contact criterion.}
For fixed $w\in B$, the total potential
\[
 \rho+v=A_0+a p(\zeta)+\lambda_0|w|^2
\]
is harmonic in $z$ on $D_\ell$. The functions
$\psi+u_0,\psi+u_1,v$ are bounded continuous
$\omega$-plurisubharmonic functions. The envelope construction and
rooftop formula in Lemmas~\ref{lem:perron} and~\ref{lem:legendre}
apply to these endpoints. Apply Theorem~\ref{thm:contact} when $n=1$,
and Remark~\ref{rem:contact-higher-dimension} when $n\ge2$, on
$D_\ell\times B$. The comparison function $v$ is smooth on this
domain, including its intersection with $E$, so the contact criterion
applies across the divisor. For the weak geodesic relative to $\omega$
joining $\psi+u_0$ and $\psi+u_1$, it gives
\[
 \varphi_t(z,w)=v(z,w)
 \qquad\text{if }m(\zeta)\le t\le m(\zeta)+h(\zeta).
\]
The strict inequalities at $\zeta_*$ give a neighborhood
$V\Subset D_\ell\times B$ of $x_*$, chosen inside $\{|z|<d/3\}$,
such that $\varphi_t=v$ near $V\times[\delta,1-\delta]$.
On $V$ we have $\psi=\varepsilon|z|^{2\beta}$, and differentiating
$\rho+v$ gives \eqref{eq:conic-main-rank}, including on $E$.
The spatial form has rank $n-1$ and null normal direction
$\partial/\partial z$. In dimension one the spatial form is zero.
\end{proof}

\begin{remark}[Tangential degeneracy]
\label{rem:conic-tangential}
When $n\ge2$, the same construction also gives tangential degeneracy
while retaining the normal cone angle of the intermediate forms.
For this variant, place the translated slit-annulus seed in the
$w_1$-coordinate in Step 1, and use
$\lambda_0(|z|^2+|w'|^2)$ in the remaining coordinates
$w'=(w_2,\ldots,w_{n-1})$. The same gluing estimates give smooth
endpoints $u_j$ and a lower function $v$ with the prescribed contacts,
whose total potential is harmonic in $w_1$ on the slit domain.

Choose $F=|s_E|_{\mathfrak h}^{2\beta}=|z|^{2\beta}$ on a neighborhood
of $\overline{\mathbb B_r}$ as above, and now set
$\widehat\psi=\varepsilon F$ without the cutoff. Add
$\widehat\psi$ to $u_0,u_1$ and $v$. For sufficiently small
$\varepsilon>0$, the background and endpoint forms are positive cone
metrics and uniformly comparable by the estimates in Step 2 with
$\chi=1$. The function $v+\widehat\psi$ is globally
$\omega$-plurisubharmonic: on $\mathbb B_r$, which contains the
support of $v$, the current $i\partial_X\bar\partial_XF$ is
nonnegative; outside that support its spatial form is the positive
cone background. Adding the same function to all three potentials
preserves their differences and every boundary contact.

Off $E$, the comparison potential $\rho+v+\widehat\psi$ is smooth
and harmonic in $w_1$. Applying
Remark~\ref{rem:contact-higher-dimension} on each transverse slice
off $E$ gives $\widehat\varphi_t=v+\widehat\psi$ near
$(V\setminus E)\times[\delta,1-\delta]$ for a neighborhood $V$ of
$x_*$. Uniqueness of bounded plurisubharmonic extension gives the
same equality across $E$. On $V\setminus E$,
\[
\begin{aligned}
 \omega+i\partial_X\bar\partial_X\widehat\varphi_t
 ={}&\bigl(\lambda_0+\varepsilon\beta^2|z|^{2\beta-2}\bigr)
       i\,dz\wedge d\bar z\\
 &+\lambda_0 i\sum_{\alpha=2}^{n-1}dw_\alpha\wedge d\bar w_\alpha.
\end{aligned}
\]
Thus the null direction is $\partial/\partial w_1$, the normal cone
angle is unchanged, and $\widehat\varphi_t-\widehat\psi=v$ is smooth
across $E$ on this region.
\end{remark}
\subsection{Poincar\'e singularities}
\label{sec:poincare-boundary}

We construct a complete Poincar\'e reference form in the prescribed
class and a fixed pair of endpoint potentials whose weak geodesic loses
rank in the normal direction on open sets approaching $E$.

\begin{proof}[Proof of Theorem~\ref{thm:poincare-boundary}]
\emph{Step 1. Choice of the reference metric.}
We first choose a complete Poincar\'e reference form in the prescribed
class with product structure near $x_*\in E$.

Choose coordinates $(z,w)$ centered at $x_*\in E$, with $E=\{z=0\}$,
and a local potential $\rho_0$ of a smooth representative $\omega_0$
of the prescribed class such that
$\rho_0(z,w)=|z|^2+|w|^2+O(|(z,w)|^3)$.
Taking $\mathfrak h$ with $|s_E|_{\mathfrak h}^2=|z|^2$ in this chart,
for fixed $a>0$ and sufficiently large $\lambda$ the initial global form
\begin{equation}\label{eq:cusp-reference}
\omega_{P,0}=\omega_0-a i\partial_X\bar\partial_X
       \log\bigl(\lambda-\log|s_E|_{\mathfrak h}^2\bigr)
\end{equation}
is a complete Poincar\'e type K\"ahler form on $X=\overline X\setminus E$
\cite[Lemma 1.1 and Proposition 1.2]{AuvrayPoincare}.
Its local potential in this chart is
\[
\rho_0(z,w)-a\log\bigl(\lambda-\log|z|^2\bigr).
\]

To construct the endpoint potentials, we hope that the global
reference form has the following local potential near $x_*$:
\begin{equation}\label{eq:cusp-local-reference}
\rho(s,w)=\rho_0(0,w)-a\log s
=|w|^2-a\log s+O(|w|^3),
\end{equation}
where $s=-\log|z|+\lambda/2$ and the remainder depends only on $w$.
We use logarithmic coordinates $\zeta=s-i\theta=-\log z+\lambda/2$,
with $\theta=\arg z$ modulo $2\pi$, for the planar construction below.
Subtracting the smooth function
\[
Q(z,w)=\rho_0(z,w)-\rho_0(0,w)
\]
from the local potential of $\omega_{P,0}$ gives
\eqref{eq:cusp-local-reference} up to an additive constant. Shrink
the chart so that the complex Hessian of $\rho_0(0,w)$ is at least
$\tfrac12 I$. The resulting local form is therefore
\[
i\partial_w\bar\partial_w\rho_0(0,w)
+\frac{a\,i\,dz\wedge d\bar z}{4|z|^2s^2}>0.
\]

To obtain a global Poincar\'e type K\"ahler form, we localize this
change in both coordinate directions. Fix a tangential radius $r_w>0$
and a smooth cutoff $0\le\chi_w\le1$, equal to one for $|w|\le r_w$
and supported in $|w|<2r_w$. Choose a smooth cutoff
$0\le\chi_z\le1$ on $\mathbb C$, equal to one for $|\xi|\le1$ and
zero for $|\xi|\ge2$. Let $r_0>0$ be a small normal radius such that
$\{|z|\le2r_0,\ |w|\le2r_w\}$ lies inside the chart, and define
\[
\chi(z,w)=\chi_z(z/r_0)\chi_w(w),\qquad
Q_{r_0}=\chi Q,\qquad
\widehat\omega=\omega_{P,0}
 -i\partial_X\bar\partial_XQ_{r_0}.
\]
Extend $Q_{r_0}$ by zero outside the chart. It is smooth on
$\overline X$. Both cutoffs equal one on
$\{|z|<r_0,\ |w|<r_w\}$, where $\widehat\omega$ has potential
\eqref{eq:cusp-local-reference}; $Q_{r_0}=0$ if $|z|\ge2r_0$
or $|w|\ge2r_w$. The normal cutoff varies only on
$r_0<|z|<2r_0$, with derivatives of order $k$ bounded by
$C_k r_0^{-k}$. The tangential cutoff and its derivatives are fixed
as $r_0$ decreases. In dimension one, omit $w$ and take $\chi_w=1$.

We check positivity on the whole support. Since $Q(0,w)=0$,
smoothness gives
\[
D_w^jQ=O(|z|)\quad(0\le j\le2),\qquad
Q_z,\ Q_{z\bar z},\ Q_{z\bar w_\alpha}=O(1),
\]
where $D_w^j$ denotes any real tangential derivative of order $j$.
The product rule yields
\begin{align*}
\partial_X\bar\partial_XQ_{r_0}
=\chi\,\partial_X\bar\partial_XQ
 +Q\,\partial_X\bar\partial_X\chi+\partial_X\chi\wedge\bar\partial_XQ
 +\partial_XQ\wedge\bar\partial_X\chi.
\end{align*}
Hence, on $0<|z|\le2r_0$,
\[
(Q_{r_0})_{z\bar z}=O(r_0^{-1}),\qquad
(Q_{r_0})_{z\bar w_\alpha}=O(1),\qquad
(Q_{r_0})_{w_\alpha\bar w_\gamma}=O(|z|),
\]
with constants independent of $r_0$. On the support,
$|z|s\le Cr_0(1+|\log r_0|)$, so these estimates give
\[
\bigl|i\partial_X\bar\partial_XQ_{r_0}\bigr|_{\omega_{P,0}}
\le C\left(\frac{|z|^2s^2}{r_0}+|z|s+|z|\right)
\le Cr_0(1+|\log r_0|^2).
\]
For sufficiently small $r_0$, the last bound is less than $1/2$.
Since $Q_{r_0}=0$ outside the support, globally
\[
\tfrac12\omega_{P,0}\le\widehat\omega
\le\tfrac32\omega_{P,0}.
\]
Hence, $\widehat\omega$ is positive and complete, and remains in the prescribed class.

\emph{Step 2. Construction of the local endpoint potentials.}
On the product neighborhood from Step~1, we construct total endpoint
potentials $q_1,q_2$ and a common lower total potential $p$.
We arrange the boundary contacts $q_j=p$ and make the corresponding
potentials $u_1,u_2$ relative to $\widehat\omega$ decay rapidly with
all Poincar\'e derivatives.

View $(s,\theta)$ as the coordinates $(x,y)$ in
Subsection~\ref{sec:noncompact}, and relabel the endpoint potentials
in~\eqref{eq:nc:endpoints}, constructed using
Lemma~\ref{lem:nc:zeros}, as $u_1^{\rm pl},u_2^{\rm pl}$. The functions
$2\theta^2+u_j^{\rm pl}$ are nonnegative, have uniformly positive
Laplacians, and have the boundary zero sets in~\eqref{eq:nc:contacts}.
The functions $u_j^{\rm pl}$ are Schwartz; on every fixed angular strip,
\eqref{eq:nc:tail-formulas} gives, for sufficiently large $s$,
$u_1^{\rm pl}=4\theta e^{-s}+2e^{-2s}$, and
$u_2^{\rm pl}=-4\theta e^{-s}+2e^{-2s}$.
We first restrict their perturbations to a fixed angular strip.
Choose $0\le\chi_Y\le1$, equal to one for $|\theta|\le Y$ and zero
for $|\theta|\ge2Y$, with $|\chi_Y^{(k)}|\le C_kY^{-k}$.
For $Y$ larger than the height of the horn,
\[
\theta^2+\tfrac12\chi_Yu_j^{\rm pl}
=(1-\chi_Y)\theta^2
 +\tfrac12\chi_Y(2\theta^2+u_j^{\rm pl})\ge0,
\]
and
\begin{align*}
\Delta\bigl(\theta^2+\tfrac12\chi_Yu_j^{\rm pl}\bigr)
=2(1-\chi_Y)
 +\tfrac12\chi_Y\Delta(2\theta^2+u_j^{\rm pl}) +\chi_Y'\partial_\theta u_j^{\rm pl}
 +\tfrac12\chi_Y''u_j^{\rm pl}.
\end{align*}
The first two terms have a fixed positive lower bound; the last two
are $O(Y^{-1})$. Thus large $Y$ preserves strict subharmonicity, and
the contacts are unchanged where $\chi_Y=1$.
Fix $b_0<\pi/4$ and a small $\sigma>0$ with $2\sigma Y<b_0$.
Replace $u_j^{\rm pl}(s,\theta)$ by
$\sigma^2\chi_Y(\theta/\sigma)
u_j^{\rm pl}(s/\sigma,\theta/\sigma)$ and rescale the horn by
$\sigma$. Keeping the notation $u_j^{\rm pl}$ for these functions,
we have
\begin{equation}\label{eq:cusp-planar-seed}
\theta^2+\tfrac12u_j^{\rm pl}\ge0,\qquad
0<c\le\Delta\bigl(\theta^2+\tfrac12u_j^{\rm pl}\bigr)\le C.
\end{equation}
Let $R$ be the sufficiently large cutoff radius in the proof of
Lemma~\ref{lem:nc:zeros}. By~\eqref{def:afaroutpotential} and
\eqref{eq:nc:tail-formulas}, the planar endpoints vanish for $s<-2R$
before rescaling. Thus $u_j^{\rm pl}=0$ for $s<-2\sigma R$, as well
as for $|\theta|\ge b_0$. They and all their derivatives decay exponentially as
$s\to+\infty$, by~\eqref{eq:nc:tail-formulas}. The zero sets of
$\theta^2+\tfrac12u_j^{\rm pl}$ on the horn boundary are inherited
from~\eqref{eq:nc:contacts}; see Figure~\ref{fig:noncompact-horn}.
The rescaled expression $\theta^2+\tfrac12u_j^{\rm pl}$ equals
$\sigma^2$ times its original value. Thus both its zeros
and its strict positivity on the opposite side are preserved.

Next, translate the rescaled horn to a neighborhood of $x_*$ with repect to the Poincare type model metric $\widehat{\omega}$.
Choose $S$ large enough that $s>S-2\sigma R$ implies $|z|<r_0$,
and let
\[
g_j(s,\theta)=\theta^2+\tfrac12u_j^{\rm pl}(s-S,\theta),
\qquad j=1,2.
\]
Let $D_S$ denote the translated horn, and continue to write
$L_0,L_1,\mathfrak C$ for its lower side, upper side and cap.
Translation preserves~\eqref{eq:cusp-planar-seed} and the contacts:
\[
\begin{aligned}
g_1&=0\ \text{on }L_0\cup\mathfrak C,
&g_1&>0\ \text{on }L_1,\\
g_2&=0\ \text{on }L_1\cup\mathfrak C,
&g_2&>0\ \text{on }L_0.
\end{aligned}
\]

The reference metric is not Euclidean in the normal direction,
so we need a different model for the total endpoint potentials.
We use $\log(1+g_j/s^2)$ and define, on $|\theta|<\pi$,
\begin{equation}\label{eq:cusp-harmonic-p}
p(w)=\rho_0(0,w),\qquad
q_j(s,\theta,w)=\rho_0(0,w)
 +\frac a2\log\left(1+\frac{g_j(s,\theta)}{s^2}\right).
\end{equation}
Since $g_j\ge0$, we have $q_j\ge p$, with $q_j=p$ exactly where
$g_j=0$. Hence the contacts are those prescribed above

We check positivity with $\Delta=\partial_s^2+\partial_\theta^2$.
After fixing the cutoff and rescaling, choose constants independent
of $S$ such that, on the strip,
\[
0\le g_j\le M,\qquad |\nabla g_j|\le M,\qquad
0<c\le\Delta g_j\le C,
\]
where $M\ge1$ and we may take $c\le2\le C$.
Differentiating $\log(1+g_j/s^2)=\log(s^2+g_j)-2\log s$ gives
\[
\Delta q_j=\frac a2\left[
\frac{\Delta g_j}{s^2+g_j}
-\frac{4s\partial_s g_j+|\nabla g_j|^2}{(s^2+g_j)^2}
+\frac{2g_j(3s^2+g_j)}{s^2(s^2+g_j)^2}
\right].
\]
Using $s^2+g_j\ge s^2$ and $3s^2+g_j\le3(s^2+g_j)$, we obtain
\[
\left|\frac{2s^2}{a}\Delta q_j-\Delta g_j\right|
\le\frac{4M}{s}+\frac{CM+M^2+6M}{s^2}.
\]
Increase $S$ so that the right side is at most $c/2$ whenever
$s\ge S-2\sigma R$. Then
\begin{equation}\label{eq:cusp-gap-positivity}
\frac{ac}{4s^2}\le\Delta q_j
\le\frac{a(C+c/2)}{2s^2}.
\end{equation}
Where the planar perturbation vanishes identically,
$g_j=\theta^2$ and $\Delta q_j=a/s^2$, so the same bounds hold
throughout the local collar.

These bounds give precisely the Poincar\'e growth in the normal
direction. Indeed, $\zeta=s-i\theta=\lambda/2-\log z$ gives
$d\zeta=-dz/z$, and there are no mixed normal--tangential terms, so
\[
i\partial_X\bar\partial_Xq_j
=\frac{\Delta q_j}{4|z|^2}\,i\,dz\wedge d\bar z
 +i\partial_w\bar\partial_w\rho_0(0,w).
\]
The tangential form is positive, and $s\sim-\log|z|$ near $E$.
Thus~\eqref{eq:cusp-gap-positivity} is equivalent to
\[
C^{-1}\widehat\omega
\le i\partial_X\bar\partial_Xq_j
\le C\widehat\omega.
\]

Since $g_j$ is uniformly bounded on the angular strip,
\[
q_j-\rho
=a\log s+\frac a2\log\left(1+\frac{g_j}{s^2}\right)
=a\log s+O(s^{-2})\qquad(s\to+\infty).
\]
Thus $q_j-\rho$ is unbounded. To obtain decaying relative potentials,
use the following local potential of the same reference form:
\[
\widetilde\rho
=\rho+\frac a2\log(s^2+\theta^2)
=\rho_0(0,w)+\frac a2\log\left(1+\frac{\theta^2}{s^2}\right).
\]
The added term is pluriharmonic on this chart, so
$i\partial_X\bar\partial_X\widetilde\rho=\widehat\omega$.
Let
\begin{align}
u_j=q_j-\widetilde\rho
&=\frac a2\log\left(
1+\frac{u_j^{\rm pl}(s-S,\theta)}{2(s^2+\theta^2)}\right),
\label{eq:cusp-relative-seed}\\
v=p-\widetilde\rho
&=-\frac a2\log\left(1+\frac{\theta^2}{s^2}\right).
\label{eq:cusp-H}
\end{align}
Then $u_j\ge v$, with the same boundary contacts.
Here $q_j$ and $\widetilde\rho$ are defined on the cut chart
$-\pi<\theta<\pi$. Since $u_j^{\rm pl}=0$ for $|\theta|\ge b_0$,
each $u_j$ vanishes near the seam and extends smoothly across it
on the punctured neighborhood.
Moreover, $\|u_j\|_\infty=O(S^{-2})$. In the tail,
\eqref{eq:cusp-relative-seed} and the exponential decay of
$u_j^{\rm pl}$ give exponential decay of $u_j$ and its ordinary
derivatives, up to powers of $s$. Hence, for every $k,N\ge0$,
\begin{equation}\label{eq:cusp-seed-decay}
|\nabla_{\widehat\omega}^{k}u_j|_{\widehat\omega} 
\le C_k\sum_{m=0}^{k}s^m|D^m u_j|
=O(s^{-N})\qquad(s\to\infty).
\end{equation}

\emph{Step 3. Global extension and boundary contacts.}
We now extend the endpoint potentials globally and construct a bounded
$\widehat\omega$-plurisubharmonic function $\hhat{v}$ below both endpoints.

Choose a smooth nondecreasing concave function $f$ on $[0,\pi]$, equal to $\theta$
for $0\le \theta\le b_0$ and constant near $\pi$, with $0\le f'\le1$.
The function $f(|\theta|)^2$ is smooth and periodic. Define
\[
\widetilde v(s,\theta)
=-\frac a2\log\left(1+\frac{f(|\theta|)^2}{s^2}\right).
\]
For $0<\theta<\pi$, direct differentiation yields
\[
\Delta\left(-\tfrac12\log(s^2+f^2)\right)
=\frac{(s^2-f^2)(1-f'^2)-ff''(s^2+f^2)}{(s^2+f^2)^2}\ge0
\quad\text{if }s>\max f.
\]
The formula extends across $\theta=\pi$ smoothly.
Thus $\rho+\widetilde v$ is plurisubharmonic, and $\widetilde v=v$
on the central sector. Since the support of $u_j$ is in
that sector, we have everywhere
\begin{equation}\label{eq:cusp-gap-order}
u_j\ge\widetilde v.
\end{equation}

Let
$T=\{s>s_0,\ |w|<r\}$ be a product collar in the chosen chart,
with $s_0>\max f$. Choose $0<\varepsilon<1/2$ and then $A>0$
small, with $4A<\varepsilon r^2$. Let $\alpha(s)$ be a smooth nonincreasing
cutoff, equal to one for $s\le1/4$ and zero for $s\ge1/2$.
To define $b$ on all of $X$, extend the local coordinate $s$ and let
\[
s=\tfrac12\bigl(\lambda-\log|s_E|_{\mathfrak h}^2\bigr),
\qquad b=2A\alpha(s/S).
\]
This $s$ agrees with $-\log|z|+\lambda/2$ in our chart, and $b$ is
bounded and smooth. Differentiating gives
\begin{align*}
i\partial_X\bar\partial_Xb
=\frac{2A}{S}\alpha'(s/S)\,i\partial_X\bar\partial_Xs+\frac{2A}{S^2}\alpha''(s/S)\,
 i\partial_Xs\wedge\bar\partial_Xs.
\end{align*}
On the support of these derivatives $s\sim S$. The estimates
$|\partial_Xs|_{\widehat\omega}=O(s)$ and
$|i\partial_X\bar\partial_Xs|_{\widehat\omega}=O(1)$ therefore give
$|i\partial_X\bar\partial_Xb|_{\widehat\omega}\le CA$,
with $C$ independent of large $S$. Choose $A$ so that $CA<1/2$.
Then $b$ is strictly $\widehat\omega$-plurisubharmonic, equals $2A$
for $s\le S/4$, and is zero for $s\ge S/2$.

On $T$ let
\[
\hat u_j=A-\varepsilon|w|^2+u_j,
\qquad \hat v=A-\varepsilon|w|^2+\widetilde v.
\]
Take $r<r_w$ and $s_0>\lambda/2-\log r_0$, so $T$ lies in the
product neighborhood of Step~1. By that step and
\eqref{eq:cusp-gap-positivity}, the forms
$\widehat\omega+i\partial_X\bar\partial_X\hat u_j$ are uniformly
comparable to $\widehat\omega$. Also $\hat v$ is
$\widehat\omega$-plurisubharmonic and $\hat v\le\hat u_j$.

Fix $A$ as above. Since $u_j^{\rm pl}(s-S,\theta)=0$ for
$s<S-2\sigma R$, formula~\eqref{eq:cusp-relative-seed} and
$|\log(1+x)|\le2|x|$ for $|x|\le1/2$ give, for large $S$,
\[
\|u_j\|_\infty
\le\frac{a\|u_j^{\rm pl}\|_\infty}{2(S-2\sigma R)^2}.
\]
Now enlarge $S$ so that this
bound is less than $A/4$.
Since $b\ge0$ and $4A<\varepsilon r^2$,
\[
\hat u_j-b\le\frac{5A}{4}-\varepsilon|w|^2<-A
\quad\text{for }3r/4<|w|<r.
\]
This region contains a neighborhood of $|w|=r$ for every $s>s_0$.
On $s_0<s<S/4$, we have $u_j=0$, $b=2A$, and hence
$\hat u_j-b=-A-\varepsilon|w|^2\le-A$.
Thus one fixed $0<\eta<A/8$ gives
$\hat v-b\le\hat u_j-b<-\eta$ near both boundary parts of $T$.
Fix a smaller transverse ball $B$ with $\varepsilon|w|^2<A/4$
on $\overline B$. On $\overline{D_S}\times\overline B$, we have $b=0$ and
\[
0\le-\widetilde v
\le\frac a2\log\left(1+\frac{b_0^2}{(S-2\sigma R)^2}\right)<A/4
\]
for large $S$, so $\hat v-b>A/2$ there.

Let $F_\eta$ be the regularized max defined
in~\eqref{eq:local-convex-cutoff}.
As in the proof of Lemma~\ref{lem:gluing}, use the reference
$\widehat\omega+i\partial_X\bar\partial_Xb$ and data
$\hat u_j-b,\hat v-b$ to define
\begin{equation}\label{eq:cusp-global-gluing}
\begin{aligned}
\hhat{u}_j&=
\begin{cases}
 b+F_\eta(\hat u_j-b),&\text{on }T,\\
 b,&\text{on }X\setminus T,
\end{cases}\\
\hhat{v}&=
\begin{cases}
 \max\{\hat v,b\},&\text{on }T,\\
 b,&\text{on }X\setminus T.
\end{cases}
\end{aligned}
\end{equation}
These are global potentials on $X$. On $X\setminus T$ they
equal $b=2A\alpha(s/S)$, hence $2A$ for $s\le S/4$ and $0$ for
$s\ge S/2$.

By~\eqref{eq:local-cutoff-Levi} and~\eqref{eq:cusp-seed-decay}, the
endpoints are smooth and bounded, with forms uniformly comparable to
$\widehat\omega$ and bounded Poincar\'e derivatives, and $\hhat{v}$ is
bounded, continuous and $\widehat\omega$-plurisubharmonic, with
$\hhat{v}\le\hhat{u}_j$. Near $\overline{D_S}\times\overline B$,
$\hhat{u}_j=\hat u_j$ and $\hhat{v}=\hat v$, so the contacts from
Step~2 are preserved.

\emph{Step 4. Completion of the proof.}
All ingredients are now ready to complete the proof.
On $D_S\times B$,
\begin{equation}\label{eq:cusp-selected-form}
\widehat\omega+i\partial_X\bar\partial_X\hhat{v}
=i\partial_w\bar\partial_w(\rho_0(0,w)-\varepsilon|w|^2).
\end{equation}
Let $m,h$ be the harmonic measures of $L_1,\mathfrak C$ in $D_S$.
The preserved contacts and the condition at infinity supplied
by~\eqref{eq:nc:end-measure} allow Theorem~\ref{thm:contact} and
Remark~\ref{rem:contact-higher-dimension} to be applied to
$\hhat{u}_1,\hhat{u}_2,\hhat{v}$ relative to $\widehat\omega$. Their weak
geodesic agrees with $\hhat{v}$ on the open region
\[
\{(\zeta,w,t)\in D_S\times B\times(0,1):
 m(\zeta)<t<m(\zeta)+h(\zeta)\},
\]
where its spatial form has rank $n-1$ and kernel spanned by
$\partial/\partial z$. On the real axis, $h>0$ and symmetry gives
$m=(1-h)/2$. As real $\zeta\to\infty$, the points
$(z,w)=(e^{\lambda/2-\zeta},0)$ approach $x_*$. Continuity therefore
gives open products $V_r\times I_r$ in this region, with
$V_r\subset X\cap B_r(x_*)$ and $1/2\in I_r\Subset(0,1)$, as required.

Finally, take the first endpoint metric as reference. The relative
endpoint potentials are then $0$ and $\hhat{u}_2-\hhat{u}_1$.
The mean value formula and chain rule applied
to~\eqref{eq:cusp-global-gluing} and~\eqref{eq:cusp-seed-decay}
give their rapid decay in every $\widehat\omega$-covariant derivative;
the preceding metric and derivative bounds give the same estimates
for this reference. Thus, we have $
 |\nabla_{\omega_P}^ku_j|_{\omega_P}=O(s^{-N})
$ as $(s\to\infty)$. We complete the proof.
\end{proof}

\begin{remark}[Weak geodesics need not remain of Poincar\'e type]
Theorem~\ref{thm:poincare-boundary} gives endpoints in
$\mathcal H_{PM}(X,\omega_P)$ whose weak geodesic does not remain in
this class. Indeed, on each constructed region $V_r\times I_r$,
\[
\omega_P+i\partial_X\bar\partial_Xu_t
=i\partial_w\bar\partial_w\bigl(\rho_0(0,w)-\varepsilon|w|^2\bigr),
\]
so the spatial form has zero normal component. Since $1/2\in I_r$
and these regions approach $x_*$, no $c>0$ can satisfy
\[
\omega_P+i\partial_X\bar\partial_Xu_{1/2}\ge c\,\omega_P
\]
throughout any neighborhood of $x_*$. In particular,
$u_{1/2}\notin\mathcal H_{PM}(X,\omega_P)$, even though both
endpoints belong to this class and have the rapid decay in the theorem.
\end{remark}
\appendix
\section{Uniqueness of weak geodesics on noncompact K\"ahler manifolds}\label{sec:bounded-comparison}

On a general noncompact K\"ahler manifold, uniqueness of bounded
weak solutions and their identification with the Perron envelope
are not automatic. Ko{\l}odziej--Nguyen
\cite[Corollary 1.5]{KolodziejNguyen2023} prove bounded comparison
on compact K\"ahler manifolds with boundary, and the second named author
\cite[Theorems 4.1 and 5.7]{YaoALE} proves comparison and
identification with the envelope in the ALE setting.
Here we record a sufficient condition for this identification on a
general K\"ahler manifold.

We consider bounded continuous endpoints
$u_0,u_1\in\PSH(X,\omega)\cap L^\infty(X)$ and weak geodesics
$u_t\in\PSH(X,\omega)\cap L^\infty(X)$, uniformly bounded for
$t\in[0,1]$. Their associated functions on $X\times\Sigma$ are
$\Omega$-plurisubharmonic solutions of~\eqref{eq:hcma} with the
prescribed joint endpoint limits.

\begin{lemma}\label{lem:bounded-uniqueness}
Let $(X,\omega)$ be a K\"ahler manifold with bounded continuous
endpoints $u_0,u_1\in\PSH(X,\omega)$. Assume that $X$ is
compact, or that there is a smooth function $b\le0$ satisfying
\begin{equation}\label{eq:uniqueness-barrier}
 \omega+i\partial_X\bar\partial_Xb\ge\tfrac12\omega,
 \qquad b\longrightarrow-\infty\quad\text{outside compact sets}.
\end{equation}
The latter condition means that, for every $M>0$, there is a compact
set $K_M\subset X$ such that $b(x)<-M$ for all $x\in X\setminus K_M$.
Then every bounded $\Omega$-plurisubharmonic solution of
\eqref{eq:hcma} with the prescribed joint endpoint limits equals the
envelope~\eqref{eq:perron}. 
\end{lemma}

\begin{proof}
In the compact case, let $b=0$ and omit conditions at spatial infinity.
Let $V$ be any bounded $\Omega$-plurisubharmonic solution of
\eqref{eq:hcma}, with the same joint endpoint limits as $U$.
It belongs to $\mathcal B$, so $V\le U$. 

It suffices to show $U\le V$.
For $0<\varepsilon \ll 1$,
choose $C>1+\max_j\|u_j\|_\infty$ and let
\[
 U_\varepsilon=(1-\varepsilon)U+\varepsilon b-\varepsilon C.
\]
This function is $\Omega$-plurisubharmonic. At the time boundary
$t=j$, its joint boundary value is
\[
 u_j+\varepsilon(b-C-u_j)\le u_j-\varepsilon,
\]
whereas $V$ has joint boundary value $u_j$.
Uniform boundedness of $U,V$ and~\eqref{eq:uniqueness-barrier} give
$U_\varepsilon\le V-\varepsilon$ outside a compact spatial set,
taking this set to be $X$ in the compact case. On that compact set,
the joint endpoint limits give fixed collars of both time boundary
circles on which $U_\varepsilon<V$. Consequently, any failure of
$U_\varepsilon\le V$ is contained in a relatively compact subset of
$X\times\Sigma$. Choose a smooth domain $D\Subset X\times\Sigma$
containing that subset, so that $U_\varepsilon\le V$ near
$\partial D$ and outside $D$.
For a smooth K\"ahler form $\Theta=\Omega+i\partial\bar\partial |\tau|^2$ near $\overline D$, the bounded
comparison principle \cite[Corollary 1.5]{KolodziejNguyen2023}
states that $f,g\in\PSH(D,\Theta)\cap L^\infty(D)$ satisfy
$f\ge g$ on $D$ if
\[
\liminf_{D\ni\xi\to\partial D}(f-g)(\xi)\ge0,
\qquad
(\Theta+i\partial\bar\partial f)^{n+1}
\le(\Theta+i\partial\bar\partial g)^{n+1}.
\]
Apply this with
\[
\Theta=\Omega+i\partial\bar\partial|\tau|^2,\qquad
f=V-|\tau|^2,\quad g=U_\varepsilon-|\tau|^2.
\]
Here $\Theta=\Omega+i\,d\tau\wedge d\bar\tau$ is K\"ahler,
and $f,g$ are bounded since $b$ is bounded on $\overline D$. Moreover,
\[
\begin{aligned}
\Theta+i\partial\bar\partial f
 &=\Omega+i\partial\bar\partial V\ge0,\\
\Theta+i\partial\bar\partial g
 &=\Omega+i\partial\bar\partial U_\varepsilon\ge0.
\end{aligned}
\]
Thus $f,g$ are $\Theta$-plurisubharmonic, and
$f-g=V-U_\varepsilon\ge0$ near $\partial D$.
Since $V$ solves~\eqref{eq:hcma},
\[
(\Theta+i\partial\bar\partial f)^{n+1}
=0\le(\Theta+i\partial\bar\partial g)^{n+1}.
\]
Comparison gives $U_\varepsilon\le V$ on $D$, hence everywhere.
Letting $\varepsilon\downarrow0$ pointwise gives $U\le V$,
and therefore $U=V$.
\end{proof}

\begin{remark}[Compact manifolds and cone singularities]
For smooth endpoints on a compact K\"ahler manifold, the compact
case of Lemma~\ref{lem:bounded-uniqueness} applies directly.
For the conic endpoints in Theorem~\ref{thm:conic-boundary}, we use
the smooth reference $\omega$ on the compact manifold $X$:
the endpoint potentials $\psi+u_0,\psi+u_1$ are bounded continuous
$\omega$-plurisubharmonic functions on $X$.
A bounded solution initially defined off $E$ extends as an
$\Omega$-plurisubharmonic function across $E\times\Sigma$.
Its Bedford--Taylor measure does not charge this divisor, so the
HCMA equation continues to hold after extension.
Thus any such bounded solution with the prescribed joint endpoint
limits equals the envelope.
\end{remark}

\begin{remark}[Euclidean and ALE ends]
Choose a nonnegative smooth function $f$ equal to $r^2$ outside a
compact set. The Euclidean asymptotics and compactness of the
remaining region give
\[
 i\partial_X\bar\partial_X f\le C\omega
\]
for some $C>0$. For $0<\varepsilon\le(2C)^{-1}$, let
$b=-\varepsilon f$. Then $b\le0$, $b\to-\infty$ outside compact
sets, and
\[
 \omega+i\partial_X\bar\partial_Xb
 \ge(1-\varepsilon C)\omega\ge\tfrac12\omega.
\]
Lemma~\ref{lem:bounded-uniqueness} therefore identifies every bounded
weak solution with the envelope. This applies to the Euclidean and
ALE examples in every complex dimension.
\end{remark}

\begin{remark}[Poincar\'e singularities]
On $X=\overline X\setminus E$, scale the smooth Hermitian norm so
that $|s_E|_{\mathfrak h}\le1$. If $\Theta_{\mathfrak h}$ denotes
the curvature form of $\mathcal O(E)$, then
\[
 b=\varepsilon\log|s_E|_{\mathfrak h}^2\le0,
 \qquad i\partial_X\bar\partial_Xb=-\varepsilon\Theta_{\mathfrak h}.
\]
The form $\Theta_{\mathfrak h}$ is smooth on $\overline X$, and
the constructed $\omega_P$ dominates a smooth compact K\"ahler
form. Hence $\Theta_{\mathfrak h}\le C\omega_P$ for some $C>0$,
and, for $0<\varepsilon\le(2C)^{-1}$,
\[
 \omega_P+i\partial_X\bar\partial_Xb\ge\tfrac12\omega_P,
 \qquad b\longrightarrow-\infty\quad\text{as }x\to E.
\]
Thus Lemma~\ref{lem:bounded-uniqueness} applies.
The endpoint estimates in Subsection~\ref{sec:poincare-boundary}
give Poincar\'e metrics in the class of
\cite[Definition 0.1]{AuvrayPoincare}. The estimate in
\cite[Theorem 2.1]{AuvrayPoincare}, together with the boundedness
of our endpoints, implies that the corresponding weak solution is
bounded. Lemma~\ref{lem:bounded-uniqueness} therefore identifies it
with our envelope.
\end{remark}
\providecommand{\bysame}{\leavevmode\hbox to3em{\hrulefill}\thinspace}
\providecommand{\MR}{\relax\ifhmode\unskip\space\fi MR }
\providecommand{\MRhref}[2]{%
  \href{http://www.ams.org/mathscinet-getitem?mr=#1}{#2}
}
\providecommand{\href}[2]{#2}

\end{document}